\documentclass[12pt,a4paper,twoside,english,american]{scrartcl}
\usepackage{helvet}
\usepackage{babel}
\usepackage{verbatim}
\usepackage{uni_math}
\usepackage{amsthm}
\usepackage{amsmath}
\usepackage{amssymb}
\usepackage[unicode=true,pdfusetitle,
 bookmarks=true,bookmarksnumbered=false,bookmarksopen=false,
 breaklinks=false,pdfborder={0 0 1},backref=false,colorlinks=false]
 {hyperref}

\makeatletter

\numberwithin{equation}{section}
\numberwithin{figure}{section}
 \theoremstyle{definition}
 \newtheorem*{defn*}{\protect\definitionname}
\theoremstyle{plain}

  \theoremstyle{remark}
  \theoremstyle{plain}
  
  \theoremstyle{plain}
  
  \theoremstyle{plain}
  
  \theoremstyle{definition}

\usepackage[T1]{fontenc}    

\usepackage{a4wide}        
\usepackage{tu-preprint}
\usepackage{amsmath}
\usepackage{enumerate}

\usepackage{euscript}
\usepackage{mathtools}
\allowdisplaybreaks[1] 

\newenvironment{keywords}{ \noindent\footnotesize\textbf{Keywords and phrases:}}{}

\newenvironment{class}{\noindent\footnotesize\textbf{Mathematics subject classification 2010:}}{}

\usepackage{color,tu-preprint}

\DeclareMathAccent{\Circ}{\mathalpha}{operators}{"17}

\renewcommand{\Re}{\operatorname{\mathfrak{Re}}}

\renewcommand*{\epsilon}{\varepsilon}
\renewcommand*{\theta}{\vartheta}
\renewcommand*{\rho}{\varrho}

\makeatother

  \addto\captionsamerican{\renewcommand{\corollaryname}{Corollary}}
  \addto\captionsamerican{\renewcommand{\definitionname}{Definition}}
  \addto\captionsamerican{\renewcommand{\lemmaname}{Lemma}}
  \addto\captionsamerican{\renewcommand{\propositionname}{Proposition}}
  \addto\captionsamerican{}
  \addto\captionsamerican{\renewcommand{\theoremname}{Theorem}}
  \addto\captionsenglish{\renewcommand{\corollaryname}{Corollary}}
  \addto\captionsenglish{\renewcommand{\definitionname}{Definition}}
  \addto\captionsenglish{\renewcommand{\lemmaname}{Lemma}}
  \addto\captionsenglish{\renewcommand{\propositionname}{Proposition}}
  \addto\captionsenglish{}
  \addto\captionsenglish{\renewcommand{\theoremname}{Theorem}}
  \providecommand{\corollaryname}{Corollary}
  \providecommand{\definitionname}{Definition}
  \providecommand{\lemmaname}{Lemma}
  \providecommand{\propositionname}{Proposition}
  
\providecommand{\theoremname}{Theorem}

\begin{document}
\selectlanguage{english}%
\institut{Institut f\"ur Analysis}

\preprintnumber{MATH-AN-09-2013}

\preprinttitle{Continuous dependence on the coefficients for a class of non-autonomous evolutionary equations}

\author{Marcus Waurick} 

\makepreprinttitlepage

\selectlanguage{american}%
\setcounter{section}{-1}

\date{}

\title{Continuous dependence on the coefficients for a class of non-autonomous evolutionary equations}

\author{Marcus Waurick\\
Institut f\"ur Analysis, Fachrichtung Mathematik\\
Technische Universit\"at Dresden\\
Germany\\
marcus.waurick@tu-dresden.de\\
}
\maketitle
\begin{abstract} The continuous dependence of solutions to certain equations on the coefficients is addressed. The class of equations under consideration has only recently be shown to be well-posed. We give criteria that guarantee that convergence of the coefficients in the weak operator topology implies weak convergence of the respective solutions. We discuss three examples: A homogenization problem for a Kelvin-Voigt model for elasticity, the discussion of continuous dependence of the coefficients for acoustic waves with impedance type boundary conditions and a singular perturbation problem for a mixed type equation. By means of counter examples we show optimality of the results obtained.
\end{abstract}
\begin{class} 35F45,  35M10,  35Q99,  46N20,  47N20, 74Q15 
\end{class}

\newpage

\tableofcontents{}


\cleardoublestandardpage

\section{Introduction}\label{sec:intro}

In this article we discuss the continuous dependence of solutions to evolutionary equations on the coefficients.
The usual method of choice to discuss issues of evolution equations is a semi-group approach. In fact many evolutionary problems in mathematical physics can be described by the abstract Cauchy problem
\[
   u'=Au\quad u(0)=u_0
\]
with $A$ being a generator of a strongly continuous semi-group in a certain Banach space $X$, $u_0\in X$. 
 Thus, in the semi-group language, we are led to consider
\[
   u_n'=A_n u_n\quad u_n(0)=u_0
\]
for a suitable sequence of generators $(A_n)_n$ in a common Banach space $X$. Now, we address whether the sequence $(u_n)_n$ of solutions converges in a particular sense. If the sequence of solutions converge to some $u$, we further ask, whether there exists an operator $A$ such that the following holds
\[
   u'=Au\quad u(0)=u_0.
\]
Within the semi-group perspective there are several issues to be taken care of: Variable domains of the generators $A_n$, non-reflexivity of the space in which the solutions $(u_n)_n$ are obtained and the generator property for $A$. To illsutrate the latter, we discuss a simple example with bounded generators: Take a bounded measurable function $a\colon \R\to \R$ and consider the sequence $(u_n)_n$ of solutions to the equation
\[
    \frac{\!\!\dd}{\!\!\dd t}u_n(t,x)+a(nx)u_n(t,x)=0\quad u(0,x)=u_0(x) \quad ((t,x)\in (0,\infty)\times \R),
\]
for some given $u_0\in L^2(\R)$. The Cauchy problem can be formulated in the state space $X=L^2(\R)$. It can be shown that, assuming for instance the peridiodicity of $a$, that the limit equation is \emph{not} of the type discussed above. Indeed, the resulting equation is of integro-differential type, see for instance \cite[Chapter 23]{TarIntro}. In particular, we cannot expect the limit equation to be of the form of the abstract Cauchy problem described above. Hence, the semi-group perspective to this kind of equation may not be advanced. The very reason for this shortcoming is that the convergence of $(a(n\cdot))_n$ is too weak, \cite{Waurick}. In fact it can be shown that $(a(n\cdot))_n$ converges in the weak star topology of $L^\infty(\R)$ to the integral mean over the period, or, equivalently, the sequence of associated multiplication operators in $L^2(\R)$ converges in the weak operator topology to the identity times the integral mean over the period of $a$. We refer to \cite{Eisner2010}, where subtleties with regards to the Trotter product formula and the weak operator topology are highlighted. Due to the non-closedness of abstract Cauchy problems with regards to the convergences under consideration, we cannot use semi-group theory here.  

To the best of the author's knowledge, besides the author's work (\cite{Waurick,Waurick2011,Waurick2012,Waurick2012a,Waurick2013a}), there are very few studies (if any) of continuous dependence on the coefficients of a general problem class under the weak operator topology. However, there are some results for particular equations and/or with stronger topologies the coefficients are considered in: In \cite{Yaman2011} a particular non-linear equation is considered and the continuous dependence of the solution on some scalar factors is addressed. Similarly, in \cite{cCelebi2006,Franchi2003,Tu2007,Payne1999,Liu2009} the so-called Brinkman-Forchheimer equation is discussed with regards to continuous dependence on some bounded functions under the sup-norm. The local sup-norm has been considered in \cite{Gianni2009}, where the continuous dependence on the (non-linear) constitutive relations for particular equations of fluid flow in porous media is discussed. A weak topology for the coefficients is considered in \cite{Kim2009}. However, the partial differential equations considered are of a specific form and the underlying spatial domain is the real line. Dealing with time-dependent coefficients in a boundary value problem of parabolic type, the author of \cite{Penning1991} shows continuous dependence of the associated evolution families on the coefficients. In \cite{Penning1991}, the coefficients are certain functions considered with the $C^1$-norm. The author of \cite{Tudor1976} studies the continuous dependence of diffusion processes under the $C^0$-norm of the coefficients. Also with regards to strong topologies, the authors of \cite{Kunze2011,Kunze2012} studied continuous dependence results for a class of stochastic partial differential equations.

  We also refer to \cite{TarIntro,CioDon,BenLiPap}, where the continuous dependence of the coefficients has been addressed in the particular situation of homogenization problems. See also the references in \cite{Waurick2013a}. Due to the specific structure of the problem semi-group theory could be applied for a homogenization problem for thermo-elasticity, \cite{Francfort1983}. 

As indicated above the main observation for discussing homogenization problems is that the coefficients might only converge in a rather weak topology. A possible choice modeling this is the weak operator topology, \cite{Waurick2012a,Waurick}. Thus, motivated by the problems in homogenization theory, we investigate the continuous dependence of solutions of evolutionary problems on the coefficients, where the latter are endowed with the weak operator topology. Aiming at an abstract result and having sketched the drawbacks of semi-group theory in this line of problems, we need to consider a different class of evolutionary equations. We focus on a certain class of integro-differential-algebraic partial differential equations. Recently, a well-posedness result could be obtained for this class, \cite{Waurick2013e}. Moreover, generalizations of the results in \cite{Waurick2012a,Waurick2012,Waurick2013a} need the development of other techniques. 

The class of equations under consideration is roughly described as follows.
Consider
\begin{equation}\label{eq:0}
  (\s Mu)' + \s A u=f,
\end{equation}
where $\s M$ is a bounded linear operator acting in space-time, $(\s Mu)'$ denotes the time-derivative of $\s Mu$ and $\s A$ is a (unbounded, linear) maximal monotone operator (see e.g.~\cite{Showalter1997}) in space-time, which is invariant under time-translations, $f$ is a given forcing term and $u$ is to be determined. The underlying Hilbert space setting will be described in Section \ref{sec:prel}. Though \eqref{eq:0} seems to be an evolution equation in any case, it is possible to choose $\s M$ in the way that \eqref{eq:0} does not contain any time-derivative at all. Indeed, in the Hilbert space framework developed below the time-derivative becomes a continuously invertible operator. Thus, as $\s M$ acts in space-time, we can choose $\s M$ as the inverse of the time-derivative times some bounded linear operator $M$, such that \eqref{eq:0} amounts to be $Mu+\s Au=f$. In view of the latter observation and in order not to exclude the algebraic type equation $Mu+\s Au=f$ we are led to consider \eqref{eq:0} with no initial data. However, imposing sufficient regularity for the initial conditions, one can formulate initial value problems equivalently into problems of the type \eqref{eq:0}, see e.g.~\cite[Section 6.2.5]{Picard}.

There are many standard equations from mathematical physics fitting in the abstract form described by \eqref{eq:0}. These are, for instance, the heat equation (\cite[Section 6.3.1]{Picard}, \cite[Theorem 4.5]{Waurick2012a}), the wave equation (\cite[Section 3]{Picard2012}, \cite[Section 4.2]{Trostorff2013}), Poisson's equation (see Section \ref{sec:ex}), the equations for elasticity (\cite[Section 4.2]{Trostorff2013}) or Maxwell's equations (\cite[Section 4.1]{PicPhy}, \cite[Section 5]{Waurick2012}). Coupled phenomena such as the equation for thermo-elasticity (\cite[Section 6.3.2]{Picard}, \cite[Theorem 4.10]{Waurick2012a}) or the equations for thermo-piezo-electro-magnetism (\cite[Section 6.3.3]{Picard}), or equations with fractional derivatives like subdiffusion or superdiffusion problems (\cite[Section 4]{Waurick2012a}, \cite[Section 4]{R.Picard2012}) can be dealt with in the general framework of \eqref{eq:0}. We note here that the operator $\s M$ in \eqref{eq:0} needs not to be time-translation invariant. Thus, the coefficients may not only contain memory terms, but they may also explicitly depend on time, see \cite[Section 3]{RainerPicard2013}. In Section \ref{sec:mixed}, we thoroughly discuss an equation of mixed type, that may involve operators $\s M$ depending on both the temporal and the spatial variable(s) or being operators of convolution type.

In order to have an idea of the form of the operators $\s M$ and $\s A$, we give three more concrete examples. All these three examples are considered in a three-dimensional spatial domain $\Omega$. Written in block operator matrix form with certain source term $f$, a first order formulation of the heat equation reads as
\[
   \left(\frac{\!\! \dd }{\!\! \dd t}\begin{pmatrix}
     1 &0 \\ 0 & 0 
   \end{pmatrix} + \begin{pmatrix}
     0 &0 \\ 0 & \kappa(t,x)^{-1} 
   \end{pmatrix} + \begin{pmatrix}
     0 & \dive \\ \grad & 0 
   \end{pmatrix}\right) \begin{pmatrix} \theta(t,x) \\ q(t,x) \end{pmatrix} = \begin{pmatrix} f(t,x) \\ 0 \end{pmatrix}, 
\]
where $\theta$ is the temperature, $q$ is the heat flux and $\kappa$ is the conductivity matrix, which is assumed to be continuously invertible. Note that the second line of the system is Fourier's law. Hence, in this case
\[
   \s M = \begin{pmatrix}
     1 &0 \\ 0 & 0 
   \end{pmatrix} + \left(\frac{\!\! \dd }{\!\! \dd t}\right)^{-1} \begin{pmatrix}
     0 &0 \\ 0 & \kappa^{-1} 
   \end{pmatrix}
\]
and 
\[
 \s A = \begin{pmatrix}
     0 & \dive \\ \grad & 0 
   \end{pmatrix},
\]
which is skew-selfadjoint if suitable boundary conditions are to be imposed.

Similarly, we find for the wave equation
\[
   \left(\frac{\!\! \dd }{\!\! \dd t}\begin{pmatrix}
     1 &0 \\ 0 & \kappa(t,x)^{-1} 
   \end{pmatrix} + \begin{pmatrix}
     0 & \dive \\ \grad & 0 
   \end{pmatrix}\right) \begin{pmatrix} u(t,x)\\ v(t,x) \end{pmatrix} = \begin{pmatrix} f(t,x) \\ 0 \end{pmatrix}, 
\]
for some suitable coefficient matrix $\kappa$, the relations
\[
  \s M = \begin{pmatrix}
     1 &0 \\ 0 & \kappa^{-1} 
   \end{pmatrix}
\]
and 
\[
 \s A = \begin{pmatrix}
     0 & \dive \\ \grad & 0 
   \end{pmatrix}.
\]
Maxwell's equations read as
\[
    \left(\frac{\!\! \dd }{\!\! \dd t}\begin{pmatrix}
     \eps(t,x) &0 \\ 0 & \mu(t,x) 
   \end{pmatrix} + \begin{pmatrix}
     \sigma(t,x) &0 \\ 0 & 0 
   \end{pmatrix} + \begin{pmatrix}
     0 & \curl \\ -\curl & 0 
   \end{pmatrix}\right) \begin{pmatrix} E(t,x) \\ H(t,x) \end{pmatrix} = \begin{pmatrix} J(t,x) \\ 0 \end{pmatrix}, 
\]
where $\eps$, $\mu$ and $\sigma$ are the material coefficients electric permittivity, magnetic permeability and the electric conductivity, respectively. $J$ is a given source term and $(E,H)$ is the electro-magnetic field. We have
\[
 \s M =\begin{pmatrix}
     \eps &0 \\ 0 & \mu 
   \end{pmatrix} + \left(\frac{\!\! \dd }{\!\! \dd t}\right)^{-1}\begin{pmatrix}
     \sigma &0 \\ 0 & 0 
   \end{pmatrix}
\]
and 
\[
 \s A = \begin{pmatrix}
     0 & \curl \\ -\curl & 0 
   \end{pmatrix},
\]
which is skew-selfadjoint for instance under the electric boundary condition. Well-posedness conditions for the above equations are suitable strict positive definiteness conditions for $\eps,\mu$ and $\kappa$. Moreover, the derivative with respect to time needs to be uniformly bounded. The precise conditions can be found in \cite[Condition (2.3)]{RainerPicard2013}.

Now, we turn to discuss the main contribution, Theorem \ref{thm:basic_conv_thm}, of the present article. Take a sequence of bounded linear operators $(\s M_n)_{n\in\N}$ in space-time converging in the weak operator topology $\tau_{\textnormal{w}}$ to some bounded linear operator $\s M$. The $\s M_n$'s are assumed to satisfy suitable conditions (see Theorem \ref{thm:Solution_theory} or \cite{Waurick2013e}) such that the respective equations as in (\ref{eq:0}) are well-posed in the sense that the (closure of the) operator $u\mapsto (\s Mu)' + \s A u$ is continuously invertible in space-time. Let $f$ be a given right-hand side. For $n\in\N$ let $u_n$ solve
\begin{equation}\label{eq:n_eq}
   (\s M_nu_n)' + \s A u_n=f.
\end{equation}
The main result now states that if the sequence of the commutator of the $\s M_n$'s with time-differentiation is a bounded sequence of bounded linear operators\footnote{If $\s M_n$ is given by multiplicaton by some function $\kappa_n$ depending on both the temporal and spatial variables, the commutator with time-differentiation is given by the operator of multiplying with the derivative of $\kappa_n$ with respect to time. Thus, the boundedness of the sequence of commutators under consideration is warranted if, for instance, $(\kappa_n)_n$ is a $C^1$-bounded sequence considered as the sequence of mappings \[\left(\R\ni t\mapsto \kappa_n(t,\cdot)\in L^\infty\right)_n.\] We also refer to Section \ref{sec:mixed} for more specific examples.} and if the resolvent of $\s A$ satisfies a certain compactness condition, we have $u_n\rightharpoonup u$, i.e., $(u_n)_n$ weakly converges to $u$, where $u$ satisfies
\[
    (\s Mu)' + \s A u=f.
\]
 It should be noted that the operator $\s A$ may only be skew-selfadjoint. In particular, the equations under consideration may not have maximal regularity. Moreover, the freedom in the choice of the sequence $(\s M_n)_n$ also allows for the treatment of differential-algebraic equations, which have applications in Control Theory, \cite{PicTroWau2012b,PicTroWau2012}. Since we only assume convergence in the weak operator topology for the operator sequence, the result particularly applies to norm-convergent sequences or sequences converging in the strong operator topology (see also Section \ref{sec:ex}). However, for the latter two cases the results are certainly not optimal. For the case of convergence in the weak operator topology, we give two examples (Examples \ref{Ex: Count} and \ref{ex:com}) that the assumptions in our main theorem cannot be dropped.

In order to proceed in equation \eqref{eq:n_eq} to the limit as $n\to\infty$, the main difficulty to overcome is to find conditions such that $(\s M_nu_n)_n$ converges to the product of the limits. This is where a compactness condition for the resolvent of $\s A$ comes into play. With this it is then possible to apply the compact embedding theorem of Aubin-Lions (see Theorem \ref{thm:Aubin-Lions} below) in order to gain a slightly better convergence of (a subsequence of) the $u_n$'s.

As it will be demonstrated in Section \ref{sec:ex}, the results have applications to homogenization theory. In a different situation, where certain time-translation invariant operators were treated, the latter has also been observed and exemplified in \cite{Waurick2012a,Waurick2012}.

Before building up the Hilbert space setting mentioned above in Section \ref{sec:prel}, we discuss a motivating example in Section \ref{sec:mixed}. Section \ref{sec:main_result} is devoted to state and briefly discuss the main result of the paper, which will be applied in Section \ref{sec:ex} to a homogenization problem in visco-elasticity, a wave equation with impedance type boundary conditions and a singular perturbation problem. The concluding section is devoted to the proof of Theorem \ref{thm:basic_conv_thm}. Any Hilbert space treated here is a complex Hilbert space. 

\section{A problem of mixed type}\label{sec:mixed}

In order to illustrate the main contribution of the present paper a bit further, we consider a partial differential equation of mixed type on a one-dimensional spatial domain. Writing $\partial_0$ for time differentiation and denoting by $\partial_1$ the distributional derivative in $L^2(0,1)$ with maximal domain and by $\interior\partial_1$ the distributional derivative in $L^2(0,1)$ with homogeneous Dirichlet boundary conditions, we treat the following system written in block operator matrix form:
\begin{multline}\label{eq:mixed_type}
   \left(\partial_0\begin{pmatrix}
              \1_{[0,\frac{1}{4}]\cup[\frac{1}{2},\frac{3}{4}]}(x) & 0 \\ 
                  0 & \1_{[0,\frac{1}{4}]\cup[\frac{3}{4},1]}(x) 
             \end{pmatrix}+ \begin{pmatrix}
              \1_{[\frac{1}{4},\frac{1}{2}]\cup[\frac{3}{4},1]}(x) & 0 \\ 
                  0 & \1_{[\frac{1}{4},\frac{1}{2}]\cup[\frac{1}{2},\frac{3}{4}]}(x) 
             \end{pmatrix}\right.\\ \left.+\begin{pmatrix} 0 & \partial_1 \\ \interior\partial_1 & 0 \end{pmatrix}\right)\begin{pmatrix} u(t,x) \\ v(t,x) \end{pmatrix} = \begin{pmatrix} f(t,x) \\ g(t,x) \end{pmatrix}\quad ((t,x)\in \R\times (0,1)), 
\end{multline}
where $\1_I$ denotes the characteristic function of some set $I$ and $f,g$ are thought of being given. In the framework sketched in the introduction, we find that $\s M$ is (use that $\partial_0$ can be established as a continuously invertible operator) given by
\[
          \s M =\begin{pmatrix}
              \1_{[0,\frac{1}{4}]\cup[\frac{1}{2},\frac{3}{4}]}(m) & 0 \\ 
                  0 & \1_{[0,\frac{1}{4}]\cup[\frac{3}{4},1]}(m) 
             \end{pmatrix}+ \partial_0^{-1}\begin{pmatrix}
              \1_{[\frac{1}{4},\frac{1}{2}]\cup[\frac{3}{4},1]}(m) & 0 \\ 
                  0 & \1_{[\frac{1}{4},\frac{1}{2}]\cup[\frac{1}{2},\frac{3}{4}]}(m) 
             \end{pmatrix},
\]
  where the argument $m$ is a reminder of multiplication-by-argument and means that the operator  $\1_{[0,\frac{1}{4}]\cup[\frac{1}{2},\frac{3}{4}]}(m)$ is to be interpreted as the multiplication operator in $L^2(0,1)$ associated with the bounded function $\1_{[0,\frac{1}{4}]\cup[\frac{1}{2},\frac{3}{4}]}$\footnote{More generally, for a function $\psi\colon\Omega \to \mathbb{C}$, $\Omega$ a set, we denote the multiplication operator induced by $\psi$ defined on $X$-valued functions for some vector space $X$ by $\psi(m)$, i.e., for $f\colon \Omega \to X$ we have $\psi(m)f\coloneqq (\omega\mapsto \psi(\omega)f(\omega))$. If $\Omega=\R$ and we want to stress that $\Omega$ has the interpretation of a 'time'-variable, we also write $\psi(m_0)$.}. Taking the boundary conditions of $\partial_1$ and $\interior\partial_1$ into account, we realize that 
\[
  \s A = \begin{pmatrix} 0 & \partial_1 \\ \interior\partial_1 & 0 \end{pmatrix}
\]
is skew-selfadjoint and, thus, maximal monotone. For further applications, it is crucial to realize that $\s A$ considered as an operator in $L^2(0,1)$ has compact resolvent. Note that the system describes a mixed type equation. Indeed, we find that the system reduces to 
\begin{itemize}
 \item $\partial_0^2u-\partial_1\interior\partial_1 u = \partial_0 f$ on $\R\times[0,\frac{1}{4}]$ with $g=0$;
 \item $\partial_0^2v-\interior\partial_1\partial_1 v = \partial_0 g$ on $\R\times[0,\frac{1}{4}]$ with $f=0$;
 \item $u-\partial_1\interior\partial_1 u =  f$ on $\R\times[\frac{1}{4},\frac{1}{2}]$ with $g=0$;
 \item $v-\interior\partial_1\partial_1 v = g$ on $\R\times[\frac{1}{4},\frac{1}{2}]$ with $f=0$;
 \item $\partial_0u-\partial_1\interior\partial_1 u = f$ on $\R\times[\frac{1}{2},\frac{3}{4}]$ with $g=0$;
 \item $\partial_0v-\interior\partial_1\partial_1 v = g$ on $\R\times[\frac{3}{4},1]$ with $f=0$. 
\end{itemize}
Hence, the system varies between hyperbolic, elliptic and parabolic type equations either with homogeneous Dirichlet or Neumann data. Well-posedness of the system \eqref{eq:mixed_type} can be established in a space-time setting. The Hilbert space the solutions of \eqref{eq:mixed_type} are obtained is $L^2_\nu(\R;L^2(0,1))$, the space of square-integrable functions with respect to the measure $\left(e^{-2\nu(\cdot)}\lambda_\R\right)\otimes \lambda_{(0,1)}$ for some $\nu>0$. Here $\lambda_I$ is the Lebesgue measure restricted to $I\subseteqq \R$. In $L^2_\nu$ it turns out that the derivative with respect to the first variable (time-variable) is normal and continuously invertible satisfying $\Re\partial_0=\nu$. Now, with the help of Theorem \ref{thm:Solution_theory} below, well-posedness of \eqref{eq:mixed_type} is warranted. Indeed, it suffices to estimate 
\[
   \Re \langle \partial_0 \s M U,\1_{\R_{\leqq a}}(m_0)U\rangle_{L_\nu^2(\R;L^2(0,1)^2)}.
\]
from below for all $U\in D(\partial_0)$ and $a\in \R$, where $\1_{\R_{\leqq a}}(m_0)$ denotes the multiplication operator induced by $\1_{\R_{\leqq a}}$ with respect to the time-variable.
We estimate for $\begin{pmatrix} u \\ v \end{pmatrix}\in D(\partial_0)$ and $a\in \R$
{\footnotesize
\begin{multline*}
\Re \langle \partial_0\s M \begin{pmatrix} u \\ v \end{pmatrix},\1_{\R_{\leqq a}}(m_0)\begin{pmatrix} u \\ v \end{pmatrix} \rangle =\\
  \left\langle \left(\nu \begin{pmatrix}
              \1_{[0,\frac{1}{4}]\cup[\frac{1}{2},\frac{3}{4}]}(m) & 0 \\ 
                  0 & \1_{[0,\frac{1}{4}]\cup[\frac{3}{4},1]}(m) 
             \end{pmatrix}+ \begin{pmatrix}
              \1_{[\frac{1}{4},\frac{1}{2}]\cup[\frac{3}{4},1]}(m) & 0 \\ 
                  0 & \1_{[\frac{1}{4},\frac{1}{2}]\cup[\frac{1}{2},\frac{3}{4}]}(m) 
             \end{pmatrix}\right)\begin{pmatrix} u \\ v \end{pmatrix} ,\1_{\R_{\leqq a}}(m_0) \begin{pmatrix} u \\ v \end{pmatrix}\right\rangle\\
\geqq \left\langle\begin{pmatrix} u \\ v \end{pmatrix},\1_{\R_{\leqq a}}(m_0)\begin{pmatrix} u \\ v \end{pmatrix}\right\rangle 
\end{multline*}}for $\nu\geqq 1$. 

Now, instead of \eqref{eq:mixed_type}, we consider the sequence of problems
 \begin{multline}
   \left(\partial_0\begin{pmatrix}
              \1_{[0,\frac{1}{4}]\cup[\frac{1}{2},\frac{3}{4}]}(n\cdot m\!\!\! \mod 1) & 0 \\ 
                  0 & \1_{[0,\frac{1}{4}]\cup[\frac{3}{4},1]}(n\cdot m \!\!\! \mod 1) 
             \end{pmatrix}\right.\\+\left. \begin{pmatrix}
              \1_{[\frac{1}{4},\frac{1}{2}]\cup[\frac{3}{4},1]}(n\cdot m \!\!\! \mod 1) & 0 \\ 
                  0 & \1_{[\frac{1}{4},\frac{1}{2}]\cup[\frac{1}{2},\frac{3}{4}]}(n\cdot m\!\!\!  \mod 1) 
             \end{pmatrix}\right. \label{eq:M_n}\\ \left.+\begin{pmatrix} 0 & \partial_1 \\ \interior\partial_1 & 0 \end{pmatrix}\right)\begin{pmatrix} u_n \\ v_n \end{pmatrix} = \begin{pmatrix} f \\ g \end{pmatrix}
\end{multline}
for $n\in \N$, where $x\!\!\!\mod 1\coloneqq x-\lfloor x\rfloor$, $x\in \R$. With the same arguments from above well-posedness of the latter equation is warranted in the space $L^2_1(\R;L^2(0,1))$. Now, \begin{multline*}
  \begin{pmatrix}
              \1_{[0,\frac{1}{4}]\cup[\frac{1}{2},\frac{3}{4}]}(n\cdot m\!\!\! \mod 1) & 0 \\ 
                  0 & \1_{[0,\frac{1}{4}]\cup[\frac{3}{4},1]}(n\cdot m \!\!\! \mod 1) 
             \end{pmatrix}\\+ \partial_0^{-1}\left. \begin{pmatrix}
              \1_{[\frac{1}{4},\frac{1}{2}]\cup[\frac{3}{4},1]}(n\cdot m \!\!\! \mod 1) & 0 \\ 
                  0 & \1_{[\frac{1}{4},\frac{1}{2}]\cup[\frac{1}{2},\frac{3}{4}]}(n\cdot m \!\!\!  \mod 1) 
             \end{pmatrix}\right.\\
\to \begin{pmatrix}
              \frac{1}{2} & 0 \\ 
                  0 & \frac{1}{2} 
             \end{pmatrix}+\partial_0^{-1} \begin{pmatrix}
              \frac{1}{2} & 0 \\ 
                  0 & \frac{1}{2} 
             \end{pmatrix}
\end{multline*}
in the weak operator topology due to peridiocity. 
Our main theorem asserts that the sequence $\begin{pmatrix} u_n \\ v_n \end{pmatrix}_n$ weakly converges to the solution $\begin{pmatrix} u \\ v \end{pmatrix}$ of the problem 
\[
   \left(\partial_0\begin{pmatrix}
              \frac{1}{2} & 0 \\ 
                  0 & \frac{1}{2} 
             \end{pmatrix}+ \begin{pmatrix}
              \frac{1}{2} & 0 \\ 
                  0 & \frac{1}{2} 
             \end{pmatrix}+\begin{pmatrix} 0 & \partial_1 \\ \interior\partial_1 & 0 \end{pmatrix}\right)\begin{pmatrix} u \\ v \end{pmatrix} = \begin{pmatrix} f \\ g \end{pmatrix}.
\]
It is interesting to note that the latter system does not coincide with any of the equations discussed above.

Our main convergence theorem deals with coefficients $\s M$ that live in space-time. Going a step further instead of treating \eqref{eq:M_n}, we let $(\kappa_n)_n$ in $W_1^1(\R)$ be a convergent sequence of weakly differentiable $L^1(\R)$-functions with limit $\kappa$ and support on the positive reals. Then it is easy to see that the associated convolution operators $(\kappa_n*)_n$ converge in $L(L^2(\R_{\geqq 0}))$ to $\kappa*$. Moreover, using Young's inequality, we deduce that
\[
   \Abs{\kappa_n*}_{L(L^2_{\nu}(\R))}, \Abs{\kappa_n'*}_{L(L^2_{\nu}(\R))} \to 0 \quad(\nu\to\infty) 
\]
uniformly in $n$. Thus, the strict positive definiteness of 
\[
   \partial_0 (1+\kappa*)\left(\begin{pmatrix}
              \1_{[0,\frac{1}{4}]\cup[\frac{1}{2},\frac{3}{4}]}(m) & 0 \\ 
                  0 & \1_{[0,\frac{1}{4}]\cup[\frac{3}{4},1]}(m) 
             \end{pmatrix}+ \partial_0^{-1}\begin{pmatrix}
              \1_{[\frac{1}{4},\frac{1}{2}]\cup[\frac{3}{4},1]}(m) & 0 \\ 
                  0 & \1_{[\frac{1}{4},\frac{1}{2}]\cup[\frac{1}{2},\frac{3}{4}]}(m) 
             \end{pmatrix}\right)
\] in the truncated form from above follows from the respective inequality for 
\[
 \partial_0\begin{pmatrix}
              \1_{[0,\frac{1}{4}]\cup[\frac{1}{2},\frac{3}{4}]}(m) & 0 \\ 
                  0 & \1_{[0,\frac{1}{4}]\cup[\frac{3}{4},1]}(m) 
             \end{pmatrix}+ \begin{pmatrix}
              \1_{[\frac{1}{4},\frac{1}{2}]\cup[\frac{3}{4},1]}(m) & 0 \\ 
                  0 & \1_{[\frac{1}{4},\frac{1}{2}]\cup[\frac{1}{2},\frac{3}{4}]}(m) 
             \end{pmatrix}.
\]
Now, the product of a sequence converging in the weak operator topology and a sequence converging in the norm topology converges in the weak operator topology. Hence, the solutions of    
\begin{multline*}
   \left(\partial_0\left(1+\kappa_n*\right)\left(\begin{pmatrix}
              \1_{[0,\frac{1}{4}]\cup[\frac{1}{2},\frac{3}{4}]}(n\cdot m\!\!\! \mod 1) & 0 \\ 
                  0 & \1_{[0,\frac{1}{4}]\cup[\frac{3}{4},1]}(n\cdot m \!\!\! \mod 1) 
             \end{pmatrix}\right.\right.\\+\left.\left. \partial_0^{-1} \begin{pmatrix}
              \1_{[\frac{1}{4},\frac{1}{2}]\cup[\frac{3}{4},1]}(n\cdot m \!\!\! \mod 1) & 0 \\ 
                  0 & \1_{[\frac{1}{4},\frac{1}{2}]\cup[\frac{1}{2},\frac{3}{4}]}(n\cdot m \!\!\!  \mod 1) 
             \end{pmatrix}\right.\right) \\ \left.+\begin{pmatrix} 0 & \partial_1 \\ \interior\partial_1 & 0 \end{pmatrix}\right)\begin{pmatrix} u_n \\ v_n \end{pmatrix} = \begin{pmatrix} f \\ g \end{pmatrix}
\end{multline*}
 converge weakly to the solution of
\[
   \left(\partial_0\left(1+\kappa*\right)\left(\begin{pmatrix}
              \frac{1}{2} & 0 \\ 
                  0 & \frac{1}{2} 
             \end{pmatrix}+ \partial_0^{-1}\begin{pmatrix}
              \frac{1}{2} & 0 \\ 
                  0 & \frac{1}{2} 
             \end{pmatrix}\right)+\begin{pmatrix} 0 & \partial_1 \\ \interior\partial_1 & 0 \end{pmatrix}\right)\begin{pmatrix} u \\ v \end{pmatrix} = \begin{pmatrix} f \\ g \end{pmatrix}.
\]
The latter considerations dealt with time-translation invariant coefficients. We shall also treat another example, where time-translation invariance is not warranted. For this take a sequence of Lipschitz continuous functions $(N_n \colon \mathbb{R}\to \R)_n$ with uniformly bounded Lipschitz semi-norm and such that $(N_n)_n$ converges pointwise almost everywhere to some function $N\colon\R\to\R$. Moreover, assume that there exists $c>0$ such that $\frac{1}{c}\geqq N_n\geqq c$ for all $n\in\N$. Then, by Lebesgue's dominated convergence theorem $N_n(m_0)\to N(m_0)$ in the strong operator topology, where we anticipated that $N_n(m_0)$ acts as a multiplication operator with respect to the temporal variable. The strict monotonicity in the above truncated sense of 
   \begin{multline*}\partial_0\left(N_n(m_0)\begin{pmatrix}
              \1_{[0,\frac{1}{4}]\cup[\frac{1}{2},\frac{3}{4}]}(n\cdot m\!\!\! \mod 1) & 0 \\ 
                  0 & \1_{[0,\frac{1}{4}]\cup[\frac{3}{4},1]}(n\cdot m \!\!\! \mod 1) 
             \end{pmatrix}\right.\\+\left.\partial_0^{-1} \begin{pmatrix}
              \1_{[\frac{1}{4},\frac{1}{2}]\cup[\frac{3}{4},1]}(n\cdot m \!\!\! \mod 1) & 0 \\ 
                  0 & \1_{[\frac{1}{4},\frac{1}{2}]\cup[\frac{1}{2},\frac{3}{4}]}(n\cdot m \!\!\!  \mod 1) 
             \end{pmatrix}\right)
\end{multline*}
is easily seen using integration by parts, see e.g.~\cite[Lemma 2.6]{RainerPicard2013}. Our main convergence theorem now yields that the solutions of  
\begin{multline*}
   \left(\partial_0\left(N_n(m_0)\begin{pmatrix}
              \1_{[0,\frac{1}{4}]\cup[\frac{1}{2},\frac{3}{4}]}(n\cdot m\!\!\! \mod 1) & 0 \\ 
                  0 & \1_{[0,\frac{1}{4}]\cup[\frac{3}{4},1]}(n\cdot m \!\!\! \mod 1) 
             \end{pmatrix}\right.\right.\\+\left.\left. \partial_0^{-1} \begin{pmatrix}
              \1_{[\frac{1}{4},\frac{1}{2}]\cup[\frac{3}{4},1]}(n\cdot m \!\!\! \mod 1) & 0 \\ 
                  0 & \1_{[\frac{1}{4},\frac{1}{2}]\cup[\frac{1}{2},\frac{3}{4}]}(n\cdot m \!\!\! \mod 1) 
             \end{pmatrix}\right.\right) \\ \left.+\begin{pmatrix} 0 & \partial_1 \\ \interior\partial_1 & 0 \end{pmatrix}\right)\begin{pmatrix} u_n \\ v_n \end{pmatrix} = \begin{pmatrix} f \\ g \end{pmatrix}
\end{multline*}
 converge weakly to the solution of
\[
   \left(\partial_0\left(N(m_0)\begin{pmatrix}
              \frac{1}{2} & 0 \\ 
                  0 & \frac{1}{2} 
             \end{pmatrix}+ \partial_0^{-1}\begin{pmatrix}
              \frac{1}{2} & 0 \\ 
                  0 & \frac{1}{2} 
             \end{pmatrix}\right)+\begin{pmatrix} 0 & \partial_1 \\ \interior\partial_1 & 0 \end{pmatrix}\right)\begin{pmatrix} u \\ v \end{pmatrix} = \begin{pmatrix} f \\ g \end{pmatrix}.
\]
The next section deals with the precise framework the solution theory is to be established.
\section{Preliminaries}\label{sec:prel}

We summarize some findings of \cite{Picard,Waurick2013d,Waurick2013e}. In the whole section let $H$ be a Hilbert space. We introduce the time-derivative operator $\partial_0$ as an operator in $L_\nu^2(\R;H)\coloneqq L^2(\R,\exp(-2\nu(\cdot))\lambda;H)$, $\lambda$ denoting the $1$-dimensional Lebesgue measure, for some $\nu>0$ as follows
\[
   \partial_0 \colon H_{\nu,1}(\R;H)\subseteqq L_\nu^2(\R;H)\to L_\nu^2(\R;H),\phi\mapsto\phi',
\]
where $H_{\nu,1}(\R;H)$ is the space of weakly differentiable $L_\nu^2(\R;H)$-functions with weak derivative also lying in the exponentially weighted $L^2$-space. For the scalar product in the latter space we occasionally write  $\langle\cdot,\cdot \rangle_\nu$. One can show that $\partial_0$ is one-to-one and that for $f\in L_\nu^2(\R;H)$ we have for all $t\in \R$ the Bochner-integral representation
\[
  \partial_0^{-1}f(t)=\int_{-\infty}^t f(\tau)\dd \tau.
\]
The latter formula particularly implies $\Abs{\partial_0^{-1}}\leqq \frac{1}{\nu}$ and, thus, $0\in \rho(\partial_0)$,  see e.g.~\cite[Theorem 2.2 and Corollary 2.5]{Kal2013} for the elementary proofs. From the integral representation for $\partial_0^{-1}$, we also read off that $\partial_0^{-1} f$ vanishes up to some time $a\in \R$, if so does $f$. This fact may roughly be described as causality. A possible definition is the following.
\begin{Def}[\cite{Waurick2013d}] Let $M\colon D(M)\subseteqq L_\nu^2(\R;H)\to L_\nu^2(\R;H)$. We say that $M$ is \emph{causal}  if for all $R>0$, $a\in \R$, $\phi\in L_\nu^2(\R;H)$ the mapping \begin{align*}
               \left( B_{M}(0,R), \abs{ \1_{\R_{\leqq a}}(m_0)\left(\cdot-\cdot\right)}\right) &\to \left(L_\nu^2(\R;H), \abs{\langle \1_{\R_{\leqq a}}(m_0)\left(\cdot-\cdot\right),\phi\rangle}\right) \\
   f&\mapsto Mf,
\end{align*} 
is uniformly continuous, where $B_M(0,R)\coloneqq \{ f\in D(M); \abs{f}+\abs{Mf}< R\}$.
\end{Def}
\begin{rems}\label{rem:causality} (a) For closed linear operators $M$, we have shown in \cite[Theorem 1.6]{Waurick2013d} that $M$ is causal if and only if for all $a\in \R$ and $\phi \in D(M)$ the implication 
\[
 \1_{\R_{<a}}(m_0)\phi=0 \Rightarrow \1_{\R_{<a}}(m_0)M\phi=0
\]
 holds. The latter, in turn, is equivalent to 
\[
   \1_{\R_{<a}}(m_0)M\1_{\R_{<a}}(m_0)=\1_{\R_{<a}}(m_0)M\quad (a\in \R)
\]
provided that $\1_{\R_{<a}}(m_0)[D(M)]\subseteqq D(M)$ for all $a\in \R$.

(b) Assume that $M\colon D(M)\subseteqq L_\nu^2(\R;H)\to L_\nu^2(\R;H)$ is continuous and for all $a\in\R$ the set 
\[
   D(M\1_{\R_{\leqq a}}(m_0))\cap D(M) \subseteqq L_\nu^2(\R;H)
\]
is dense\footnote{The latter happens to be the case if, for instance, $D(M)\supseteqq C_{\infty,c}(\R;H)$, the space of indefinitely differentiable functions with compact support.}. If for all $a\in \R$ we have
\[
   \1_{\R_{<a}}(m_0)M\1_{\R_{<a}}(m_0)=\1_{\R_{<a}}(m_0)M
\]
on $D(M\1_{\R_{<a}}(m_0))\cap D(M)$, then both $M$ and $\overline{M}$ are causal. Indeed, by continuity, the latter equality implies that
\begin{align*}
  \1_{\R_{<a}}(m_0)\overline{M}\1_{\R_{<a}}(m_0)&=\overline{\1_{\R_{<a}}(m_0)M\1_{\R_{<a}}(m_0)}\\
                                                &=\overline{\1_{\R_{<a}}(m_0)M}=\1_{\R_{<a}}(m_0)\overline{M} \quad (a\in \R).
\end{align*}
 Hence, by (a) $\overline{M}$ is causal, implying the causality for $M$.
\end{rems}

\begin{rems}\label{rem:mat_law_function}
  A prototype of causal operators are particular functions of $\partial_0^{-1}$.\footnote{In \cite[Section 4]{Waurick2013a}, \cite{Waurick}, \cite{R.Picard2012} examples for this kind of operators are given. Some of these are convolutions with suitable $L^1$-functions or the time-shift.}  Though being of independent interest, we need this class of operators to properly formulate the examples in Section \ref{sec:ex}. We use the explicit spectral theorem for $\partial_0^{-1}$ given by the \emph{Fourier-Laplace transformation} $\s L_\nu$. Here $\s L_\nu$ is the unitary transformation from $L_\nu^2(\R)$ onto $L^2(\R)$ such that
\[
   \s L_\nu f = \left(x\mapsto \frac{1}{\sqrt{2\pi}}\int_\R e^{-ixy-\nu y}f(y)\dd y\right)
\]
 for continuous functions $f$ with compact support. Then, one can show that
\[
   \partial_0^{-1}= \s L_\nu^* \frac{1}{im_0+\nu}\s L_\nu,
\]
where $\frac{1}{im_0+\nu}\phi(x)\coloneqq \frac{1}{ix+\nu}\phi(x)$ for $\phi\in L^2(\R)$, $x\in \R$. Now, any $M$ belonging to the Hardy-space $\s H^\infty(B(r,r))$ of bounded and analytic functions $B(r,r)\to \mathbb{C}$ for some $r>\frac{1}{2\nu}$ leads to a causal, time-translation invariant operator $M\left(\partial_0^{-1}\right)$ in the way that
\[
   M\left(\partial_0^{-1}\right) \coloneqq \s L_\nu^* M\left(\frac{1}{im_0+\nu}\right)\s L_\nu.
\]
We endow $\s H^\infty(B(r,r))$ with the supremum norm.
 Moreover, note that the definitions made can readily be extended to the vector-valued case, i.e., if $\s H^\infty(B(r,r);L(H))$ denotes the Hardy space of bounded analytic functions with values in the space of bounded linear operators, we can define for $M\in \s H^\infty(B(r,r);L(H))$ the operator
\begin{equation}\label{eq:mat_law_fn}
   M\left(\partial_0^{-1}\right) \coloneqq \s L_\nu^* M\left(\frac{1}{im_0+\nu}\right)\s L_\nu,
\end{equation}
 acting in the Hilbert space $L_\nu^2(\R;H)$, where we re-used $\s L_\nu$ to denote the extension of the scalar-valued Fourier-Laplace transformation to the $H$-valued one. Thus, \eqref{eq:mat_law_fn} should be read in the strong sense. With the help of the Paley-Wiener theorem, it is possible to show causality for $M(\partial_0^{-1})$, see e.g.~\cite{PicPhy}.
\end{rems}

In \cite{Waurick2013e}, we have shown the following well-posedness result, which comprises a large class of linear partial integro-differential-algebraic equations of mathematical physics as it has been demonstrated in \cite{Picard,RainerPicard2013} (see also Section \ref{sec:intro}). Before we state the theorem, we introduce the notion of a bounded commutator.

\begin{Def}[bounded commutator] Let $B\in L(H)$ and  let $A\colon D(A)\subseteqq H\to H$ be a densely defined linear operator. Then $B$ is said to have a \emph{bounded commutator with $A$}, if there exist $C\in L(H)$ such that $BA\subseteqq AB+C$. In the latter case, we shall write $[B,A]\coloneqq -[A,B]\coloneqq C$. A sequence $(B_n)_n$ of bounded linear operators is said to have \emph{bounded commutators with $A$}, if for all $n\in \N$ the operator $B_n$ has a bounded commutator with $A$ and the sequence $([B_n,A])_n$ is bounded. 
\end{Def}

\begin{Sa}[{\cite[Theorem 2.7 and 2.4]{Waurick2013e}}]\label{thm:Solution_theory} Let $H$ be a Hilbert space, $\s M\in L(L_\nu^2(\R;H))$. Assume that $\s M$ has a bounded commutator with $\partial_0$.
 Let $\s A\colon D(\s A)\subseteqq L_\nu^2(\R;H)\to L_\nu^2(\R;H)$ be linear, maximal monotone and such that $\partial_0 (\s A+1) = (\s A+1) \partial_0 $, i.e., $\s A$ \emph{commutes with }$\partial_0$.
Moreover, assume the positive definiteness conditions
\begin{equation}\label{ineq:thm:Solth}
   \Re \left\langle \partial_0 \s M u,\1_{\R_{\leqq a}}(m_0)u\right\rangle\geqq c\left\langle u,\1_{\R_{\leqq a}}(m_0)u \right\rangle, \quad \Re \left\langle \s A u,\1_{\R_{\leqq 0}}(m_0)u\right\rangle\geqq 0
\end{equation}
for all $u\in D(\partial_0)\cap D(\s A)$, $a\in \R$, and some $c>0$.

Then $0\in \rho\left(\overline{\partial_0 \s M+\s A}\right)$ and the operator $\left(\overline{\partial_0 \s M+\s A}\right)^{-1}$ is causal.
\end{Sa}
For the following we need to record some continuity estimates. In order to do so, we briefly recall the concept of Sobolev lattices discussed in \cite[Chapter 2]{Picard}. 
\begin{Def}[(short) Sobolev lattice] Let $C_1,C_2$ be two densely defined closed linear operators in $H$ with $0\in \rho(C_1)\cap \rho(C_2)$ and $C_1C_2=C_2C_1$. Then, for $k,\ell\in \{-1,0,1\}$ we define $H_{k,\ell}(C_1,C_2)$ as the completion of $D(C_1C_2)$ with respect to the norm $\phi\mapsto \abs{C_1^{k}C_2^{\ell}\phi}$. The family $(H_{k,\ell}(C_1,C_2))_{k,\ell\in\{-1,0,1\}}$ is called \emph{(short) Sobolev lattice}. 
\end{Def}
 
\begin{rems}\label{rem:SoLa}
  (a) We have continuous embeddings $H_{k,\ell}(C_1,C_2)\hookrightarrow H_{k',\ell'}(C_1,C_2)$ provided that $k\geqq k'$ and $\ell\geqq \ell'$. Moreover, $H_{k,\ell}(C_1,C_2)=H_{\ell,k}(C_2,C_1)$ for all $k,k',\ell,\ell'\in \{-1,0,1\}$.
 
  (b) The operators $C_1^{\pm 1}$ can be established as unitary operators from the space $H_{k,\ell}(C_1,C_2)$ into $H_{k\mp1,\ell}(C_1,C_2)$ for all $k,\ell\in \{-1,0,1\}$ such that $k\mp1\in \{-1,0,1\}$ and similarly for $C_2$.

  (c) In the special case of $C_2=1$, we write $H_{k}(C_1)\coloneqq H_{k,1}(C_1,1)$ for all $k\in\{-1,0,1\}$.

  (d) In the special case of $C_1=\partial_0$ and $C_2=1$, we write $H_{\nu,k}(\R;H)\coloneqq H_{k}(\partial_0)$ for all $k\in\{-1,0,1\}$.  
\end{rems}

\begin{rem}\label{rem:domainoftime-space}
  With the help of the Sobolev lattice construction stated, we can drop the closure bar in $\overline{\partial_0\s M+\s A}$ and compute in the Sobolev lattice associated with $(\partial_0,\s A+1)$. In order to make this more precise, we denote here the extensions of $\partial_0$ and $\s A$ to the Sobolev lattice (with the common domain $H_{0,0}(\partial_0,\s A+1)=L_\nu^2(\R;H)$) by $\partial_0^{\textnormal{e}}$ and $\s A^{\textnormal{e}}$, respectively. Now, let $u\in D(\overline{\partial_0\s M+\s A})\subseteqq L_\nu^2(\R;H)$. Then, by definition, there exists a sequence $(u_n)_n$ in $D(\partial_0\s M)\cap D(\s A)$ such that $u_n\to u$ and $v_n\coloneqq \left(\partial_0\s M+\s A\right)u_n \to \overline{\left(\partial_0\s M+\s A\right)}u\eqqcolon v$ in $L_\nu^2(\R;H)$ as $n\to\infty$. On the other hand, the continuity of $\partial_0^{\textnormal{e}}$ and $\s A^{\textnormal{e}}$ implies $\partial_0^\textnormal{e}\s M u_n \to \partial_0^\textnormal{e}\s M u$ and $\s A^\textnormal{e} u_n \to \s A^\textnormal{e} u$ in $H_{-1,0}(\partial_0,\s A+1)$ and $H_{0,-1}(\partial_0,\s A+1)$, respectively, as $n\to\infty$. From $L_\nu^2(\R;H)\hookrightarrow H_{-1,-1}(\partial_0,\s A+1)$ and  
\[
  v_n = \partial_0^\textnormal{e}\s M u_n + \s A^\textnormal{e} u_n \stackrel{n\to\infty}{\to} \partial_0^\textnormal{e}\s M u + \s A^\textnormal{e} u \in H_{-1,-1}(\partial_0,\s A+1)
\]
it follows that $v=\partial_0^\textnormal{e}\s M u + \s A^\textnormal{e} u$. Thus,
\[
   D(\overline{\partial_0\s M+\s A})\subseteqq\{ u\in L_\nu^2(\R;H); \partial_0^\textnormal{e} \s Mu +\s A^\textnormal{e}u\in L_\nu^2(\R;H)\}\eqqcolon D
\]
and $(\overline{\partial_0\s M+\s A})u=\partial_0^\textnormal{e} \s Mu +\s A^\textnormal{e}u$ for $u\in D(\overline{\partial_0\s M+\s A})$. 

On the other hand, if $u\in D$ then one can show that $(1+\eps\partial_0^{-1})u \in D(\partial_0)\cap D(\s A)\subseteqq D(\partial_0 \s M+\s A)$ for every $\eps>0$, see \cite[Lemma 4.2]{Waurick2013e} or \cite[Lemma 2.9]{RainerPicard2013}. Moreover, from the lemmas stated, it also follows that $\left(\left(\partial_0 \s M+\s A\right)(1+\eps\partial_0^{-1})u\right)_{\eps>0}$ is weakly convergent in $L_\nu^2(\R;H)$ as $\eps\to 0+$. As the strong closure of linear operators coincides with the weak closure, we deduce that $u\in D(\overline{\partial_0\s M+\s A})$.
\end{rem}

With the observations made in the latter remark, we henceforth omit the closure bar in $\overline{\partial_0\s M+\s A}$, use the continuous extensions of $\partial_0$ and $\s A$ to the Sobolev lattice, re-use the respective notations and agree that $D(\partial_0 \s M+\s A)=\{ u\in L_\nu^2(\R;H); \partial_0 \s Mu +\s Au\in L_\nu^2(\R;H)\}$.
Now, we are in the postition to state the continuity estimates:
\begin{Fo}\label{cor:cont_estim} In the situation of Theorem \ref{thm:Solution_theory}, let $f\in L_\nu^2(\R;H)$ and let $u\in L_\nu^2(\R;H)$ be the solution of $(\partial_0\s M + \s A)u=f$. Then $u\in H_{-1,1}(\partial_0,\s A+1)$ and 
\[
  \abs{u}_{H_{-1,1}(\partial_0,\s A+1)}\leqq \left(\frac{1}{\nu}+\Abs{\s M}_{L(L_\nu^2(\R;H))}\frac{1}{c}  + \frac{1}{c\nu}\right)\abs{f}_{L_\nu^2(\R;H)}
\]
\end{Fo}
\begin{proof}
  In the Sobolev lattice associated to $(\partial_0,\s A +1)$ we compute that
 \[
    (\s A+1)u = f-\partial_0 \s M u + u \in H_{-1,0}(\partial_0,\s A+1).
 \]
Thus, 
\[
   u = (\s A+1)^{-1}\left(f-\partial_0 \s M u + u \right) \in H_{-1,1}(\partial_0,\s A+1)
\]
and
\begin{align*}
   \abs{u}_{-1,1}&=\abs{\partial_0^{-1}(\s A+1) u}_{0,0}\\
                 &=\abs{\partial_0^{-1}(\s A+1) (\s A+1)^{-1}\left(f-\partial_0 \s M u + u \right)}_{0,0}\\
                 &\leqq\abs{\partial_0^{-1}f}_{0,0}+\abs{\s M u}_{0,0} + \abs{\partial_0^{-1}u}_{0,0}\\
                 &\leqq\frac{1}{\nu}\abs{f}_{0,0}+\Abs{\s M}\frac{1}{c} \abs{f}_{0,0} + \frac{1}{c\nu}\abs{f}_{0,0}\\
                 &\leqq\left(\frac{1}{\nu}+\Abs{\s M}\frac{1}{c}  + \frac{1}{c\nu}\right)\abs{f}_{0,0}\qedhere
\end{align*}
\end{proof}

\section{The basic convergence theorem}\label{sec:main_result}

We recall the concept of $G$-convergence:

\begin{Def}[{$G$-convergence, {\cite[p.\ 74]{Gcon1}, \cite{Waurick2012}}}]
 Let $H$ be a Hilbert space. Let $(A_{n}:D(A_{n})\subseteqq H\to H)_{n}$
be a sequence of continuously invertible linear operators onto $H$
and let $B:D(B)\subseteqq H\to H$ be linear and one-to-one. We say
that $(A_{n})_{n}$ \emph{$G$-converges to $B$} if $\left(A_{n}^{-1}\right)_{n}$
converges in the weak operator topology to $B^{-1}$, i.e., for all
$f\in H$ the sequence $(A_{n}^{-1}(f))_{n}$ converges weakly to
some $u$, which satisfies $u\in D(B)$ and $B(u)=f$. $B$ is called
the%
\footnote{Note that the $G$-limit is uniquely determined, cf. \cite[Proposition 4.1]{Waurick2012}.%
} \emph{$G$-limit} of $(A_{n})_{n}$ and we write $A_{n}\stackrel{G}{\longrightarrow}B$. 
\end{Def}
Our main theorem reads as follows.
\begin{Sa}\label{thm:basic_conv_thm} Let $H$ be a Hilbert space, $\nu>0$. Let $\left(\s M_n\right)_n$ be a bounded sequence in $L(L_\nu^2(\R;H))$ with bounded commutators with $\partial_0$. Moreover, let $\s A\colon D(\s A)\subseteqq L_\nu^2(\R;H)\to L_\nu^2(\R;H)$ linear and maximal monotone commuting with $\partial_0$ and assume that $\s M_n$ is causal, $n\in\N$.
Moreover, assume the positive definiteness conditions
\[
  \Re \left\langle \partial_0 \s M_n u,\1_{\R_{<a}}(m_0)u\right\rangle\geqq c\left\langle u,\1_{\R_{<a}}(m_0)u \right\rangle,\quad \left\langle \s A u,\1_{\R_{<0}}(m_0)u\right\rangle\geqq 0
\]
for all $u\in D(\partial_0)\cap D(\s A)$, $a\in \R$, $n\in\N$ and some $c>0$.

Assume that there exists a Hilbert space $K$ such that $K\hookrightarrow\hookrightarrow H$, i.e., $K$ is compactly embedded into $H$, $H_1(\s A+1)\hookrightarrow L_\nu^2(\R;K)$ and that $\left(\s M_n\right)_n$ converges in the weak operator topology to some $\s M$.

Then $\partial_0\s M+\s A$ is continuously invertible in $L_\nu^2(\R;H)$ and $\left(\partial_0\s M_n+\s A\right)\stackrel{G}{\to}\left(\partial_0\s M+\s A\right)$ as $n\to\infty$.
\end{Sa}

\begin{rems}\label{rem:stronger_conv} It should be noted that it is possible to show another continuity property. Namely, if $(f_n)_n$ in $L_\nu^2(\R;H)$ is a weakly convergent sequence with\footnote{We denote the support of a function $v\colon \R\to X$ with values in some topological vector space $X$ by $\spt v$.} $\inf_n\inf\spt f_n>-\infty$ and $(u_n)_n$ is the sequence of solutions to 
\[
   (\partial_0 \s M_n+\s A)u_n = f_n,
\]
then $(u_n)_n$ weakly converges to the solution $u$ of 
\[
  (\partial_0\s M+\s A)u = \textnormal{w-}\lim_{n\to\infty} f.
\] 
\end{rems}

In view of the well-posedness theorems \cite[Theorem 2.13]{RainerPicard2013} and \cite[Theorem 2.4]{Waurick2013e}, there is a more adapted version of Theorem \ref{thm:basic_conv_thm}:

\begin{Fo}\label{co:basic_conv_thm} Let $H$ be a Hilbert space, $\nu>0$. Let $\left(\s M_n\right)_n, \left(\s N_n\right)_n$ be  bounded sequences of causal operators in $L(L_\nu^2(\R;H))$ having bounded commutators with $\partial_0$ and $\s A\colon D(\s A)\subseteqq L_\nu^2(\R;H)\to L_\nu^2(\R;H)$ linear, maximal monotone commuting with $\partial_0$.
Assume the positive definiteness conditions
\begin{align*}
  \Re \left\langle \left(\partial_0 \s M_n+\s N_n\right)u,\1_{\R_{<a}}(m_0)u\right\rangle & \geqq c\left\langle u,\1_{\R_{<a}}(m_0)u \right\rangle \quad (a\in \R)\\
  \Re \left\langle \s A u ,\1_{\R_{<0}}(m_0) u\right\rangle &\geqq 0
\end{align*}
for all $u\in D(\partial_0)\cap D(\s A)$,  $n\in\N$ and some $c>0$.

Assume that there exists a Hilbert space $K$ such that $K\hookrightarrow\hookrightarrow H$ and $H_1(\s A+1)\hookrightarrow L_\nu^2(\R;K)$ and that $\left(\s M_n\right)_n$, $\left(\s N_n\right)_n$ converges in the weak operator topology to some $\s M$ and $\s N$, respectively.

Then $\partial_0\s M+\s N+\s A$ is continuously invertible in $L_\nu^2(\R;H)$ and $\left(\partial_0\s M_n+\s N_n+\s A\right)\stackrel{G}{\to}\left(\partial_0\s M+\s N+\s A\right)$ as $n\to\infty$.
\end{Fo}

\begin{proof}
 It suffices to verify the assumptions in Theorem \ref{thm:basic_conv_thm} on $(\s M_n)_n$ for the operator sequence $\left(\s M_n+\partial_0^{-1}\s N_n\right)_n$. This, however, is easy to see.
\end{proof}

We remark here that in order to prove well-posedness for equations of the form 
\[
   (\partial_0 \s M + \s N + \s A)u = f
\]
for suitable $\s M$, $\s N$, $\s A$ in \cite[Theorem 2.13]{RainerPicard2013} or \cite[Theorem 2.4]{Waurick2013e}, we did not need any assumptions on the commutator of $\s N$ and $\partial_0$. Thus, one might wonder, whether the boundedness for the commutators of $\left(\s N_n\right)_n$ with $\partial_0$ is needed in Corollary \ref{co:basic_conv_thm}. The next example shows that this boundedness assumption is needed to compute the limit equation in the way it is done in Corollary \ref{co:basic_conv_thm}.

\begin{Bei}[On the boundedness of {$([\s N_n,\partial_0])_n$}]\label{ex:com} Let $\nu>0$. Consider for $n\in\N$ the operator 
\[
   \sin(n m_0) \colon L_\nu^2(\R) \to L_\nu^2(\R), f\mapsto \left(\sin(n\cdot)f(\cdot)\right).
\]
Define for $n\in\N$ the operators $\s M_n=0$, $\s N_n \coloneqq  \sin(n m_0) + 2$ and $\s A \colon \mathbb{C} \to \mathbb{C}, x\mapsto x$. Then, clearly, the (uniform) positive definiteness condition is satisfied and $\s A$ has compact resolvent. For $f\in C_{\infty,c}(\R)$ consider the problem of finding $u_n\in L_\nu^2(\R)$ such that
\[
   \left(\partial_0\s M_n + \s N_n +\s A\right) u_n =f,
\]
  which is the same as to say that
\[
   \left(\left(\sin(n m_0)+2\right) + 1 \right)u_n =f.
\]
We get that $u_n =\frac{1}{\sin(nm_0)+3}f$. By periodicity of $\sin$ we get with the help of  \cite[Theorem 2.6]{CioDon},
that 
\[
  u_n \rightharpoonup \int_{-\pi}^\pi \frac{1}{\sin(t)+3}\dd t f= \frac{\pi}{\sqrt{2}}f\eqqcolon u, 
\]
as $n\to\infty$. Moreover, it is easy to see that
\[
   \s N_n \stackrel{\tau_{\textnormal{w}}}{\to} \int_{-\pi}^\pi \sin(t)+2\dd t = 4\pi\eqqcolon \s N.
\]
Thus, if the representation formulas for the limit equation remain true also in this case, we would obtain that $u$ satisfies the equation
\[
  (4\pi+2) \frac{\pi}{\sqrt{2}}f + 1 \frac{\pi}{\sqrt{2}}f = \s N u + \s A u = f,
\]
which is not true, since $4\pi+3\neq \frac{\sqrt{2}}{\pi}$. A reason for this is that we cannot deduce that the weak limit of the sequence $(\s N_n u_n)_n$ equals the product of the respective limits. Indeed, we have
\[
 \s N  u = (4\pi+2) \frac{\pi}{\sqrt{2}}f \neq \textnormal{w-}\lim_{n\to\infty} \s N_n u_n =\left(\int_{-\pi}^\pi \frac{\sin(t)+2}{\sin(t)+3}\dd t\right)f =\pi\left(2-\frac{1}{\sqrt{2}}\right)f.
\]
\end{Bei}

Though the latter example does not fit into the scheme developed above it well fits into the theory established in \cite{Waurick2013a}, where we did not need the assumptions on the sequence having bounded commutator with $\partial_0$. 

We recall \cite[Example 4.9]{Waurick2012} to show that the compactness condition on $\s A$ is also needed to compute the limit in the way it is done in \ref{co:basic_conv_thm}.

\begin{Bei}[Compactness assumption does not hold]\label{Ex: Count} Let $\nu,\eps>0$. Consider the mapping $a:\R \to\R$ given by
\[
   a(x) := \1_{[0,\frac12)}(x-k)+2\1_{[\frac12,1]}(x-k)
\]
for all $x\in  [k,k+1)$, where $k\in\Z$. Define the corresponding multiplication operator in $L_2(\R)$, i.e.\ for $\phi \in C_{\infty,c}(\R)$, $a(n\cdot m)\phi:=(x\mapsto a(nx)\phi(x))$ for $n\in\N$. Note that $a(x+k)=a(x)$ for all $x\in\R$ and $k\in\Z$.
Let $f\in L^2_{\nu}(\R; L^2(\R))$. We consider the evolutionary equation with $(\s M_n)_n\coloneqq(0)_n$, $(\s N_n)_n\coloneqq(a(n\cdot m))_n$ and $\s A=i: L^2(\R)\to L^2(\R): \phi\mapsto i\phi$. By \cite[Theorem 2.6]{CioDon}, we deduce that 
\[
   \s N_n\to \frac{3}{2}
\]
as $n\to\infty$. If the assertion of Theorem \ref{thm:basic_conv_thm} remains true in this case, then $(\s N_n+\s A)_n$ $G$-converges to $\frac32 + i$. For $n\in\N$, let $u_n\in L^2_{\nu}(\R;L^2(\R))$ be the unique solution of the equation
\begin{equation}\label{a_n}
  (\s N_n +\s A)u_n=(a(n m )+i)u_n=f.
\end{equation}
Observe that by \cite[Theorem 2.6]{CioDon}
\[
  u_n=(a(nm)+i)^{-1}f\rightharpoonup \left(\int_0^1(a(x)+i)^{-1}\dd x\right) f=:u.
\]
as $n\to\infty$. We integrate
\[
  \int_0^1(a(x)+i)^{-1}\dd x =\frac12(1+i)^{-1}+\frac12(2+i)^{-1}.
\]
Inverting the latter equation yields
\[
  \left(\int_0^1(a(x)+i)^{-1}\dd x\right)^{-1} =\left(\frac12(1+i)^{-1}+\frac12(2+i)^{-1}\right)^{-1}=\frac{18}{13}+\frac{14}{13}i.
\]
Hence, $u$ satisfies
\[
  \left(\frac32 + i\right)u=f \text{ and }\left(\frac{18}{13}+\frac{14}{13}i\right)u=f,
\]
which of course is a contradiction.
\end{Bei}

\section{Examples}\label{sec:ex}

\subsection*{A time-dependent Kelvin-Voigt model}

We discuss an Example from \cite{Abdessamad2009}, we also refer to \cite{Amosov2013}, where some convergence estimates have been established. For showing that this example fits into our abstract scheme, we introduce some operators first.
Let $\Omega\subseteqq \R^3$ be open and bounded. Denote the weak symmetrized gradient acting on square-integrable vector fields in $L^2(\Omega)^3$ with (generalized) Dirichlet boundary condition by $\interior{\Grad}$. Korn's inequality implies $D(\interior{\Grad})=H_{1,0}(\Omega)^3$. By definition, for $v\in D(\interior{\Grad})$ the mapping $\interior{\Grad}v$ is an element of $H_{\textnormal{sym}}(\Omega)$, the space of square-integrable symmetric $3\times 3$-matrices. Endowing the latter space with the inner product
\[
  (\Phi,\Psi)\mapsto \int_\Omega \trace(\Phi(x)^*\Psi(x)) \dd x,
\]
we realize that $\interior{\Grad}=-\Div$, where the latter operator is the weak row-wise divergence with maximal domain. Note that by Korn's inequality and Rellich's selection theorem, we have that $D(\interior\Grad)\hookrightarrow\hookrightarrow L^2(\Omega)$, where the first space is endowed with the graph-norm of $\interior{\Grad}$. From Poincare's inequality, we see that $\interior{\Grad}$ has closed range. Denote by $\iota_R \colon R(\interior{\Grad}) \to H_{\textnormal{sym}}(\Omega)$ the canonical injection. As a consequence, the operator $\iota_R^*$ is the orthogonal projection onto the range of $\interior{\Grad}$, see e.g.~\cite[Lemma 3.2]{R.Picard2012}. 

In order to treat the problem class properly, we need to recall some notions from \cite{Waurick2013e} and \cite{Waurick2013a}:
\begin{Def}[evolutionary mappings, {\cite[Definition 2.1]{Waurick2013a}}] Let $\nu_1>0$. For Hilbert spaces $H_{0},H_{1}$, we call a linear
mapping 
\begin{equation}
M\colon D(M)\subseteqq\bigcap_{\nu\geqq\nu_1}L_{\nu}^{2}(\mathbb{R};H_{0})\to\bigcap_{\nu\geqq\nu_{1}}L_{\nu}^{2}(\mathbb{R};H_{1})\label{eq:evolutionary_map}
\end{equation}
\emph{evolutionary (at $\nu_{1}$)} 
if $D(M)\subseteq L_{\nu}^{2}(\mathbb{R};H_0)$ is dense and $M\colon D(M)\subseteq L_\nu^2(\R;H_0)\to L_\nu^2(\R;H_1)$ is closable for all $\nu\geqq\nu_{1}$.
We say $M$ is \emph{bounded}, if, in addition, $M_\nu\coloneqq \overline{M}\in L\left(L_{\nu}^{2}(\mathbb{R};H_{0}),L_{\nu}^{2}(\mathbb{R};H_{1})\right)$ for all $\nu\geqq\nu_{1}$ such that%
\footnote{For a bounded linear operator $A$ from $L_{\nu}^{2}(\mathbb{R};H_{0})$ to
$L_{\nu}^{2}(\mathbb{R};H_{1})$ we denote its operator norm by $\left\Vert A\right\Vert _{L(L_{\nu}^{2}(\mathbb{R};H_{0}),L_{\nu}^{2}(\mathbb{R};H_{1}))}$.
If the spaces $H_{0}$ and $H_{1}$ are clear from the context, we
shortly write $\left\Vert A\right\Vert _{L(L_{\nu}^{2})}$. %
} 
\[
\limsup_{\nu\to\infty}\left\Vert M\right\Vert _{L(L_{\nu}^{2})}<\infty.
\]
We define 
\[
  L_{\text{ev},\nu_{1}}(H_{0},H_{1})\coloneqq\{M;M\text{ is as in (\ref{eq:evolutionary_map}), is evolutionary at }\nu_{1}\text{ and bounded}\}
\]
and abbreviate $L_{\text{ev},\nu_{1}}(H_{0})\coloneqq L_{\text{ev},\nu_{1}}(H_{0},H_{0})$.
We call $\mathfrak{M}\subseteqq L_{\text{ev},\nu_{1}}(H_{0},H_{1})$ \emph{bounded} if $\limsup_{\nu\to\infty}\sup_{M\in\mathfrak{M}}\|M\|_{L(L_{\nu}^{2})}<\infty.$
 A family $(M_{\iota})_{\iota\in I}$ in $L_{\text{ev},\nu_{1}}(H_{0},H_{1})$
 is called \emph{bounded }if $\{M_{\iota};\iota\in I\}$ is bounded. 
\end{Def}
In \cite{Waurick2013a,Waurick2013e}, we gave several examples for evolutionary mappings. Multiplication operators are a particular subclass of these. Moreover, the operator $\s A\coloneqq \begin{pmatrix} 0 & \Div\iota_R \\ \iota_R^*\interior{\Grad} & 0\end{pmatrix}$ (defined in space-time) is also evolutionary for every $\nu>0$ and even bounded evolutionary in $L_{\textnormal{ev},\nu}(D(\Div\iota_R)\oplus D(\iota_R^*\interior{\Grad});R(\interior{\Grad})\oplus L^2(\Omega)^3)$. Trivially, $\s A$ is causal. For bounded evolutionary mappings, we recall the following result:
 \begin{Le}[{\cite[Lemma 3.3]{Waurick2013e}}] Let $\nu\geqq\nu_1\geqq \nu_0$, $H_0,H_1$ Hilbert spaces. Let $M\in L_{\textnormal{ev},\nu_0}(H_0,H_1)$ be causal. Then $M_\nu$ and $M_{\nu_1}$ coincide on $L_{\nu_1}^2(\R;H_0)\cap L_\nu^2(\R;H_0)$. 
\end{Le} 

In view of the latter lemma, we omit the subscript in the notation of the closures for causal, evolutionary mappings for different values of $\nu$, if there is no risk of confusion.

Now, take $\nu>0$ and let $\rho\in L_{\textnormal{ev},\nu}(L^2(\Omega)^3)$, $A,B\in L_{\textnormal{ev},\nu}(L^2(\Omega)^{3\times 3})$. The model treated in \cite{Abdessamad2009}, can be written as 
\[
  \partial_0 \rho \partial_0 u - \Div B \interior\Grad \partial_0 u - \Div A \interior\Grad u = f. 
\]
Abbreviating $v\coloneqq \partial_0 u$ and using $\Div B \interior\Grad = \Div \iota_R \iota_R^* B \iota_R\iota_R^*\interior{\Grad}$ (see e.g.~\cite{Waurick2012a}), we arrive at 
\[
  \partial_0 \rho v - \Div \iota_R \left(\iota_R^* \left(B +A \partial_0^{-1}\right)\iota_R\right)\iota_R^*\interior{\Grad} v = f.
\]
Now, if $B_\nu$ is \emph{strictly positive definite (uniformly for all large $\nu$)} in the sense that 
\[
  \Re \langle B_\nu u,\1_{\R_{<a}}(m_0)u\rangle \geqq c\langle u,\1_{\R_{<a}}(m_0)u\rangle
\]
 for some $c>0$ and all sufficiently large $\nu$, $u\in D(B)$, $a\in \R$ and $\nu$ is chosen large enough, we end up with ($q\coloneqq  \left(\iota_R^* \left(B +A \partial_0^{-1}\right)\iota_R\right)\iota_R^*\interior{\Grad}v$):
\[
   \left(\partial_0 \begin{pmatrix} \rho & 0 \\ 0 & 0 \end{pmatrix} + \begin{pmatrix} 0 & 0 \\ 0 & \left(\iota_R^* \left(B +A \partial_0^{-1}\right)\iota_R\right)^{-1} \end{pmatrix} - \begin{pmatrix} 0 & \Div\iota_R \\ \iota_R^*\interior{\Grad} & 0 \end{pmatrix}\right)\begin{pmatrix} v \\ q \end{pmatrix} =\begin{pmatrix} f \\ 0 \end{pmatrix}. 
\]
Assuming that $\rho,A,B$ have bounded commutators with $\partial_0$ in the sense of Theorem \ref{thm:Solution_theory}\footnote{Note that the boundedness of the commutator of $A$ and $\partial_0$ is not needed to ensure the well-posedness of the respective equation. For general well-posedness conditions for this particular equation we refer to the concluding section in \cite{RainerPicard2013}.}. If, in addition, $\partial_0\rho_\nu$ is strictly positive definite (uniformly for all large $\nu$), then it is easy to see that the aformentioned Kelvin-Voigt model for visco-elasticity is well-posed in the sense of Theorem \ref{thm:Solution_theory}. Moreover, it is easy to see that if $A,B$ and $\rho$ are thought of as being multiplication operators, the assumption on the boundedness of the commutator follows if one assumes that the respective functions are Lipschitz continuous and almost every where strongly differentiable (with respect to the temporal variable). For the latter see \cite{RainerPicard2013,Waurick2013e}.
  Thus,
\[
 \left(\partial_0 \begin{pmatrix} \rho & 0 \\ 0 & 0 \end{pmatrix} + \begin{pmatrix} 0 & 0 \\ 0 & \left(\iota_R^*\left( B +A \partial_0^{-1}\right)\iota_R\right)^{-1} \end{pmatrix} - \begin{pmatrix} 0 & \Div\iota_R \\ \iota_R^*\interior{\Grad} & 0 \end{pmatrix}\right)
\]
is continuously invertible in the underlying Hilbert space $L_\nu^2(\R;L^2(\Omega)^3\oplus R(\interior{\Grad}))$. As well-posedness issues are not the focus of the present article, we now apply our abstract homogenization theorem:
\begin{Sa} Let $\nu>0$, $(\rho_n)_n$, $(B_n)_n$, $(A_n)_n$ be bounded sequences of causal operators in $L_{\textnormal{ev},\nu}(L^2(\Omega)^3)$, $L_{\textnormal{ev},\nu}(H_\textnormal{sym}(\Omega))$, and $L_{\textnormal{ev},\nu}(H_\textnormal{sym}(\Omega))$, respectively. Assume that the respective sequences have bounded commutators with $\partial_0$. Moreover, assume there exists $c>0$ such that
\[
   \Re \langle B_n u,\1_{\R_{\leqq a}}(m_0)u\rangle_{\nu'}\geqq c\langle \phi,\1_{\R_{\leqq a}}(m_0)\phi\rangle_{\nu'},\quad \Re \langle \partial_0 \rho_n u,\1_{\R_{\leqq a}}(m_0)u\rangle_{\nu'}\geqq c\langle u,\1_{\R_{\leqq a}}(m_0)u\rangle_{\nu'} 
\]
for all $\nu'\geqq \nu$ and $\phi\in L_\nu^2(\R;H_{\textnormal{sym}}(\Omega))$ and $u\in H_{\nu,1}(\R;L^2(\Omega)^3)$, $a\in \R$.

Then there exists a subsequence $(n_k)_k$ such that
\begin{multline*}
  \left(\partial_0 \begin{pmatrix} \rho_{n_k} & 0 \\ 0 & 0 \end{pmatrix} + \begin{pmatrix} 0 & 0 \\ 0 & \left(\iota_R^* \left(B_{n_k} +A_{n_k} \partial_0^{-1}\right)\iota_R\right)^{-1} \end{pmatrix} - \begin{pmatrix} 0 & \Div\iota_R \\  \iota_R^*\interior{\Grad} & 0 \end{pmatrix}\right) \stackrel{G}{\to}\\
\left(\partial_0 \begin{pmatrix} \rho & 0 \\ 0 & 0 \end{pmatrix} + \begin{pmatrix} 0 & 0 \\ 0 & \sum_{\ell=0}^\infty \mathcal{M}_\ell \end{pmatrix} - \begin{pmatrix} 0 & \Div\iota_R \\ \iota_R^*\interior{\Grad} & 0\end{pmatrix}\right), 
\end{multline*}
where the latter operator is continuously invertible and \[\mathcal{M}_\ell=\tau_{\textnormal{w}}\textnormal{-}\lim_{k\to\infty} \left(-\left(\iota_R^* B_{n_k} \iota_R\right)^{-1}\iota_R^*A_{n_k}\iota_R\partial_0^{-1} \right)^\ell\left(\iota_R^* B_{n_k} \iota_R\right)^{-1}\] and $\rho=\tau_\textnormal{w}\textnormal{-}\lim_{k\to\infty}\rho_{n_k}$. 
\end{Sa}
\begin{proof}
 The proof follows with a Neumann series expansion of $\left(\iota_R^* \left(B_{n_k} +A_{n_k} \partial_0^{-1}\right)\iota_R\right)^{-1}$, the fact that for a sequence $(T_n)_n$ converging in the weak operator topology in some Hilbert space $H$, we have $\Abs{\tau_{\textnormal{w}}\textnormal{-}\lim_{n\to\infty} T_n}\leqq \liminf_{n\to\infty}\Abs{T_n}$, that for a separable Hilbert space $H$ norm-bounded subsets of $L(H)$ are relatively compact and metrizable with respect to the weak operator topology, and Theorem \ref{thm:basic_conv_thm}. 
\end{proof}

\begin{rems} (a) We give some more explicit formulae for $\s M_\ell$ for particular situations:
\begin{enumerate}[(i)]
 \item In the particular case, where $A_n=0$ and $B_n$ is time-independent, i.e., for every $n\in \N$ there exists $b_n\in L(H_{\textnormal{sym}}(\Omega) )$ such that $B_n$ is the (canonical) extension of $b_n$ to $L_\nu^2(\R;H_{\textnormal{sym}}(\Omega))$, then $\s M_\ell=0\, (\ell\in \N_{>0})$ and $\s M_0 =\lim_{k\to\infty} (\iota_R^* b_{n_k}\iota_R)^{-1}$. One can show that in the special case of $b_n=d(n\cdot)$, where $d$ is a matrix of suitable size with entries in the space of $[0,1]^3$-periodic $L^\infty(\R^3)$-functions, the result coincides with the classical limit, which will be addressed in a future publication. For the classical result see e.g.~\cite{TarIntro}.  
 \item Assume that $(B_n)_n=(b_n)_n$, where $(b_n)_n$ is a bounded sequence of causal operators in $L_{\textnormal{ev},\nu}(\mathbb{C})$ such that $(b_n)_\nu \geqq c$ for all $n\in\N$ and some $c>0$ and such that the sequence of respective commutators with $\partial_0$ is bounded as well. Furthermore, assume that $(A_n)_n=(a_n)_n$, where $a_n=d(n\cdot)$ for a function $d$ as in the previous part. Now, if $b_n\to b$ strongly\footnote{Here strong convergence means that there exists $b\in L_{\textnormal{ev},\nu}(\mathbb C)$ such that for any $\nu'\geqq \nu$ we have that $(b_n)_{\nu'}\to (b)_{\nu'}$ as $n\to\infty$ in the strong operator topology.} for some $b\in L_{\textnormal{ev},\nu}(\mathbb C)$, then 
\[
   \mathcal{M}_\ell=\tau_{\textnormal{w}}\textnormal{-}\lim_{k\to\infty} \left(-\iota_R^*A_{n_k}\iota_R\right)^\ell \left(b^{-1}\partial_0^{-1} \right)^\ell b^{-1} \quad (\ell\in\N).
\]
 \item Assume that both $(A_n)_n$ and $(B_n)_n$ satisfy the structural assumption on $B_n$ as in (ii) being representable as operators only acting in time. Assume, in addition, that $B_n$ is uniformly strictly positive (as in (ii)) and that $((A_n)_n,(B_n^{-1})_n,(\partial_0^{-1})_n)$ has the \emph{product convergence property}, see \cite[Definition 5.1]{Waurick2013a}, then 
\[
  \mathcal{M}_\ell=\iota_R^*\left(\tau_{\textnormal{w}}\textnormal{-}\lim_{n\to\infty} \left(-B_n^{-1}A_{n}\partial_0^{-1}\right)^\ell \left(B_n\right)^{-1}\right)\iota_R \quad (\ell \in \N).
\] 
\end{enumerate}
(b) We shall note here that the considerations above can be done similarly for the case of $B_n=0$ and $A_n$ selfadjoint and (uniformly) strictly positive definite. The homogenization result then coincides with the classical one in the sense of part (a)(i).

(c) There is also a possibility to treat the cases (a)(i) and (b) in a unified way. The resulting formulas, however, become more involved. We refer to the concluding section in \cite{RainerPicard2013} for a unified treatment of the cases (a)(i) and (b) with regards to well-posedness issues.
\end{rems}

\subsection*{The wave equation with impedance type boundary conditions}

We recall the setting in \cite[Section 3]{Picard2012} or \cite[Section 4]{Trostorff2013}. We let $\Omega\subseteqq \R^n$ be a bounded open set such that $H_1(\Omega)\hookrightarrow\hookrightarrow L^2(\Omega)$,\footnote{There is a vast literature on compact embedding theorems for the space of weakly differentiable $L^2(\Omega)$-functions into $L^2(\Omega)$. In order to maintain such compact embedding, one has to assume some 'regularity' property of the boundary of $\Omega$, see e.g.~\cite{Adam,Wloka1982}} i.e.,  the maximal domain of the distributional gradient $\grad$ defined on $L^2(\Omega)$ endowed with the graph norm of $\grad$ is compactly embedded into $L^2(\Omega)$. Analogously let $\diverg$ be the distributional divergence on $L^2(\Omega)^n$ with maximal domain. The respective skew-adjoints will be denoted by $\interior{\dive}$ and $\interior{\grad}$, as these operators encode homogeneous Neumann and Dirichlet boundary conditions, respectively.

Formally, the equations treated in \cite[Section 3]{Picard2012} (or in \cite[Section 4]{Trostorff2013}) read as
\[
   \left(\partial_0 \s M + \begin{pmatrix}
                      0 & \dive \\ \grad & 0
                     \end{pmatrix}\right)U=F,
\]
for some given $F$ and $\s M$. We address the continuous dependence on the coefficient $\s M$. Imposing additional structure on $\s M$ and the right-hand side $F$, we may rewrite the latter system into a more common form. Indeed, if $F=(f,0)$ and $\s M=\diag(\s M_1,\s M_2)$ with respect to the block structure of $\begin{pmatrix}
                      0 & \dive \\ \grad & 0
                     \end{pmatrix}$ we obtain with $U=(u_1,u_2)$:
\[
   \partial_0 \s M_1 u_1 + \dive u_2 = f \text{ and }\partial_0 \s M_2 u_2 + \grad u_1 = 0,
\]
which leads to
\begin{equation}\label{eq:wave_imp}
   \partial_0 \s M_1 u_1 - \dive \s M_2^{-1}\partial_0^{-1} \grad u_1 = f. 
\end{equation}

Choosing an appropriate domain for $\begin{pmatrix}
                      0 & \dive \\ \grad & 0
                     \end{pmatrix}$ in \emph{space-time}, which will be done below, it is possible to show that $\partial_0^{-1}\grad u=\grad\partial_0^{-1}u$ for suitable $u$. Thus, equation \eqref{eq:wave_imp} reads
\[
   \partial_0 \s M_1 u_1 - \dive \s M_2^{-1}\grad \partial_0^{-1} u_1 = f.
\]
Substituting $u\coloneqq \partial_0^{-1}u_1$, we arrive at
\begin{equation}\label{eq:wave_truely}
    \partial_0 \s M_1\partial_0 u - \dive \s M_2^{-1}\grad u = f,
\end{equation}
which may be regarded as the wave equation in a more familiar form.

Before we address continuous dependence of the solution on the coefficients, we comment on the choice of the domain for $\begin{pmatrix} 0 & \dive \\ \grad & 0  \end{pmatrix}$.

 Let $\nu>0$. As in \cite[Section 3]{Picard2012}, we take a time-translation-invariant subspace of the maximal domain of 
\[
  \begin{pmatrix} 0 & \dive \\ \grad & 0 \end{pmatrix} \subseteqq \left(L_{\nu}^2(\R;L^2(\Omega)\oplus L^2(\Omega)^n)\right)^2,
\]
such that the respective operator satisfies the conditions imposed on $\s A$ in Theorem \ref{thm:Solution_theory}. 
For this, we let $r>\frac{1}{2\nu}$ and $a\colon B(r,r)\to L^\infty(\Omega)^n$ bounded, analytic. Similar to Remark \ref{rem:mat_law_function}, $a$ gives rise to an operator in $L(L_\nu^2(\R;L^2(\Omega),L_\nu^2(\R;L^2(\Omega)^n)$ in the way that if 
\begin{equation}\label{eq:power_series_a}
   a(z) = \sum_{k=0}^\infty a_{k,r}(m)(z-r)^k\quad (z\in B(r,r))
\end{equation}
is the power series expression for $a$ in $r$ for suitable $L^\infty(\Omega)^n$-elements $a_{k,r}$, we define
\begin{align*}
  a(\partial_0^{-1})\phi &\coloneqq \sum_{k=0}^\infty a_{k,r}(m)(\partial_0^{-1}-r)^k\phi \\
                         &\coloneqq \sum_{k=0}^\infty ((x,t)\mapsto a_{k,r}(x)\left((\partial_0^{-1}-r)^k\phi\right)(t,x) 
\end{align*}
  for $\phi\in C_{\infty,c}(\R\times \Omega)$.

Throughout, we assume the following smoothness conditions on the coefficients in \eqref{eq:power_series_a}:
The mappings
\begin{align*}
   z\mapsto & \dive a (z)\coloneqq \sum_{k=0}^\infty \left(\dive a_{k,r}\right)(m)(z-r)^k\\
   z\mapsto & \curl a (z)\coloneqq \sum_{k=0}^\infty \left(\curl a_{k,r}\right)(m)(z-r)^k
\end{align*}
are bounded, analytic with $\dive a_{k,r}$ and $\curl a_{k,r}$ being measurable and bounded functions (the latter condition of course only in the case $n=3$).

Now, we are in the position to define the domain metioned above:\footnote{We shall note here that in \cite[p. 541]{Picard2012} the condition $a(\partial_0^{-1})\phi-\psi\in L_\nu^2(\R;D(\interior{\dive}))$ is replaced by  $a(\partial_0^{-1})\phi-\partial_0^{-1}\psi\in L_\nu^2(\R;D(\interior{\dive}))$.}
\[
  D(\s A) \coloneqq  \left\{ (\phi,\psi)\in D\left( \left(\begin{smallmatrix} 0&\dive \\ \grad & 0\end{smallmatrix}\right)\right) ; a(\partial_0^{-1})\phi-\psi \in L_\nu^2(\R;D(\interior{\dive}))\right\},
\]
with $\s A=\left(\begin{smallmatrix} 0&\dive \\ \grad & 0\end{smallmatrix}\right)$ on $D(\s A)$. 
We shall note here that the boundary conditions introduced include the Robin boundary conditions or boundary conditions with temporal convolutions at the boundary, cf.~\cite[p. 542]{Picard2012}. By the definition, we see that $\s A$ is time-translation invariant. Henceforth, we will also impose the sign-constraint (\cite[formula (3.3)]{Picard2012}) on $a(\partial_0^{-1})$:
\begin{equation}\label{eq:AisMaxMon}
   \Re \int_{-\infty}^0 \left(\langle \grad p, a(\partial_0^{-1})p\rangle (t)+\langle p,\dive a(\partial_0^{-1})p\rangle (t)\right)e^{-2\nu t}\dd t\geqq 0
\end{equation}
for all $p\in L^2_{\nu}(\R;D(\grad))$.
We have
\begin{Sa}\label{thm:a_maxmon} The operator $\s A$ is maximal monotone. Moreover, we have
\[
   \Re \langle \s Au,\1_{\R_{<0}}(m_0)u\rangle \geqq 0 \quad (u\in D(\s A)).
\]
If, in addition, we have that 
\begin{equation}\label{eq:a_sksa}
   \Re \int_\R \left(\langle \grad p, a(\partial_0^{-1})p\rangle (t)+\langle p,\dive a(\partial_0^{-1})p\rangle (t)\right)e^{-2\nu t}\dd t = 0
\end{equation}
for all $p\in L_{\nu}^2(\R;D(\grad))$, then $\s A$ is skew-selfadjoint.
\end{Sa}

Before the proof we record the following fact communicated by Sascha Trostorff.
\begin{Prop}\label{prop:Amm} Let $H$ be a Hilbert space, $A\colon D(A)\subseteqq H\to H$ linear. Assume that both $A$ and $-A$ are maximal monotone. Then $A$ is skew-selfadjoint. 
\end{Prop}
\begin{proof}
 From $\Re \langle Au,u\rangle\geqq 0$ and $ \Re \langle- Au,u\rangle\geqq 0$, it follows that $\Re \langle Au,u\rangle=0$ for all $u\in D(A)$. Thus, by polarization, $-A\subseteqq A^*$. The maximal monotonicity of $A$ implies the (maximal) monotonicity for $A^*$. The maximality of $-A$ yields $-A=A^*$.
\end{proof}

\begin{proof}[Proof of Theorem \ref{thm:a_maxmon}]
 \cite[Proposition 3.2]{Picard2012} shows the inequality stated and the closedness of $\s A$. Time-translation invariance together with \cite[Proposition 3.3]{Picard2012}, which for $u\in D(A^*)$ asserts that
\[
   \Re \langle\s A^*u,\1_{\R_{<0}}(m_0)u\rangle \geqq 0,
\]
 yields
\[
  \Re \langle\s Au,u\rangle,\Re \langle\s A^*v,v\rangle \geqq 0\ (u\in D(\s A),v\in D(\s A^*)).
\]
The latter together with the closedness of $\s A$ implies the maximal monotonicity for $\s A$.

Now, assume the validity of \eqref{eq:a_sksa}. Then the above reasoning shows that both $\s A$ and $-\s A$ are maximal monotone. The assertion follows from Proposition \ref{prop:Amm}.
\end{proof}

In view of Theorem \ref{thm:Solution_theory}, we also need the following result:
\begin{Prop}\label{prop:comm_with} Let $H$ be a Hilbert space, $\nu>0$ and $\s B\colon D(\s B)\subseteqq L_\nu^2(\R;H)\to L_\nu^2(\R;H)$ densely defined, closed, linear with $0\in \rho(\s B)$. Assume that $\tau_h\s B=\s B\tau_h$ for all $h\in \R$ on $D(\s B)$, where $\tau_h \in L(L_\nu^2(\R;H))$ with $\tau_h f\coloneqq f(\cdot+h)$. Then $\partial_0^{-1}(\s B)^{-1} = (\s B)^{-1}\partial_0^{-1}$.
\end{Prop}
\begin{proof}
  For $h\in\R\setminus\{0\}$ and $u\in L_\nu^2(\R;H)$, we have
\[
 \frac{1}{h}(\tau_h-1)(\s B)^{-1}\partial_0^{-1}u= (\s B)^{-1}\frac{1}{h}(\tau_h-1)\partial_0^{-1}u.
\]
Now, since $\partial_0^{-1}u\in D(\partial_0)$ and $(\s B)^{-1}$ is continuous the right-hand-side converges to $(\s B)^{-1}u$ as $h\to 0$. Thus, the left-hand side is bounded, weak compactness of $L_\nu^2$ now implies that the left-hand side converges weakly, the limit equals $\partial_0 (\s B)^{-1}\partial_0^{-1}u$. The assertion follows.
\end{proof}
Now, from $\partial_0^{-1}(\s A+1)^{-1} = (\s A+1)^{-1}\partial_0^{-1}$ and $0\in \rho(\partial_0)\cap\rho(\s A+1)$ it follows that $\partial_0(\s A+1)=(\s A+1)\partial_0$, see e.g.~\cite[p. 56]{Picard}, \cite[Lemma 1.1.1]{Waurick2011}.

In order to show a continuous dependence result on the coefficients, we need to warrant the compactness condition for the operator $\s A$ in Theorem \ref{thm:basic_conv_thm}. For higher dimensions, the nullspace of the operator $\s A$ discussed in this section is infinite-dimensional. Thus, if we want to apply Theorem \ref{thm:basic_conv_thm}, we have to consider the reduced operator $\iota_N^*\s A\iota_N$, where $\iota_N \colon N(\s A)^\bot \to L_{\nu}^2(\R;L^2(\Omega)\oplus L^2(\Omega)^n)$ is the canonical embedding from the orthogonal complement of the nullspace of $\s A$ into $L_{\nu}^2(\R;L^2(\Omega)\oplus L^2(\Omega)^n)$. The latter procedure of course is not needed, if we restrict ourselves to the one-dimensional case:
\begin{Sa}\label{thm:imp_1dim}  Let $\nu>0$. Assume that $\Omega$ is a bounded, open interval, and let $\left(\s M_k\right)_k$ be a sequence of causal operators in $L_\nu^2(\R;L^2(\Omega)\oplus L^2(\Omega))$ converging in the weak operator topology such that the sequence has bounded commutators with $\partial_0$. If, in addition, there exists $c>0$ such that
\[
   \Re \langle \partial_0\s M_n u,\1_{\R_{<a}}(m_0) u\rangle \geqq c\langle u,\1_{\R_{<a}}(m_0)u\rangle \quad (n\in \N,a\in \R,u\in D(\partial_0))
\]
then
\[
   \partial_0 \s M_n + \s A \stackrel{G}{\to} \partial_0 \s M + \s A
\]
in $L_\nu^2(\R;L^2(\Omega)^2)$.
\end{Sa}
\begin{proof}
 For the proof note that the Hilbert space $D(\partial_1)\oplus D(\partial_1)=D(\grad)\oplus D(\dive)=H_1(\Omega)^2$ is compactly embedded into $L_2(\Omega)^2$. Moreover, the validity of the conditions in Theorem \ref{thm:basic_conv_thm} are easily checked with the help of Theorem \ref{thm:a_maxmon} and Proposition \ref{prop:comm_with}. Thus, Theorem \ref{thm:basic_conv_thm} applies.
\end{proof}

\begin{rems}
 With the second order formulation of equation \eqref{eq:wave_truely}, we consider
\begin{equation*}
    \partial_0 \s M_{1,n}\partial_0 u_n - \dive \s M_{2,n}^{-1}\grad u_n = f.
\end{equation*}
for $(\s M_{1,n})_n$, $(\s M_{2,n})_n$ being such that $\s M_n\coloneqq \diag(\s M_{1,n},\s M_{2,n})$ satisfies the assumptions of Theorem \ref{thm:imp_1dim}. It follows that $\s M_{j,n}\stackrel{\tau_{\textnormal{w}}}{\to} \s M_{j}$ for some $\s M_{j}$, $j\in \{1,2\}$. The limit equation would then be the following
\[
    \partial_0 \s M_{1}\partial_0 u - \dive \s M_{2}^{-1}\grad u = f.
\]
We note here that at first one computes the limit of $(\s M_{2,n})_n$ and after that one inverts the limit to get the latter equation. In the classical terms, i.e., under certain structural and peridicity assumptions,  $\s M_2^{-1}$ is the harmonic mean of the $\s M_{2,n}^{-1}$'s. 
\end{rems}

%

Next, we discuss whether the compactness property assumed in Theorem \ref{thm:basic_conv_thm} for $\s A$ holds in the case of dimension $n=3$, which will be assumed in the remainder of this section. Recall that our strategy relies on considering the reduced operator $\iota_N^*\s A\iota_N$. We state a first important consequence:
\begin{Prop}\label{prop:red_ismax}
  The operator $\iota_N^*\s A\iota_N$ is maximal monotone. If $\s A$ is skew-selfadjoint, then so is $\iota_N^*\s A\iota_N$.
\end{Prop}
\begin{proof}
 It is plain that the operator is monotone. Thus, by Minty's theorem, it suffices to show that $1+\iota_N^*\s A\iota_N$ is onto. For this let $y\in N(\s A)^\bot$. By the maximal monotonicity of $\s A$, there exists $x\in D(\s A)$ such that $x+\s Ax = y$. We multiply the latter equality by $\iota_N^*$, which gives $\iota_N^*x+\iota_N^*\s A x =\iota_N^*y=y$. Decomposing $x=x_1+x_2$ for some $x_1\in N(\s A)^\bot$ and $x_2\in N(\s A)$, we get that $\iota_N^*x=x_1$ and $\s A x=\s A(x_1+x_2)=\s Ax_1 = \s A \iota_N\iota_N^*x$. Hence, $\iota_N^*x$ is the desired element in the domain of $1+\iota_N^*\s A \iota_N$ mapped to $y$.
The last assertion of the proposition, follows from Proposition \ref{prop:Amm}.
\end{proof}
As a next step, we need to verify that $\iota_N^*\s A\iota_N$ satisfies the assumptions in our main homogenization theorem. For this, however, we need to impose additional regularity of the boundary of $\Omega$. With addditional effort, these regularity requirements can certainly be relaxed. Since we are only interested in providing a class of examples rich enough, we do not follow the way of presenting a streamlined version of a particular compactness result.
\begin{Sa}\label{thm:definitionofK} Assume, in addition, that $\Omega$ is of class $C_5$. Then, we have, that
\[
   H_{1}(1+\iota_N^*\s A\iota_N)\hookrightarrow L_\nu^2(\R;H_1(\Omega)^4),
\]
\end{Sa}

Before we go into the proof of the theorem, we state the main ingredient: Gaffney's inequality. For the latter recall the operator $\curl$ being the distributional $\curl$ defined on $L^2(\Omega)^3$ with values in $L^2(\Omega)^3$ with maximal domain. We also use the canonical extension of $\curl$ to space-time and re-use the notation. It will become clear from the context which operator is used.
 
\begin{Sa}[{Gaffney's inequality, see e.g.~\cite[below Theorem 8.6, p.~157]{Leis1986}}]\label{thm:gaffney} Let $\Omega$ belong to the class $C_5$. Then there exists $c>0$ such that for all $u\in D(\curl)\cap D(\interior{\dive})$ we have
\[
   \abs{u}_{H_1(\Omega)}\leqq c\left(\abs{u}_{L^2(\Omega)} +\abs{\dive u}_{L^2(\Omega)}+\abs{\curl u}_{L^2(\Omega)}\right).
\] 
\end{Sa}

\begin{proof}[Proof of Theorem \ref{thm:definitionofK}] At first, we observe that
\[
 \begin{pmatrix}
    0 & \dive|_{C_{\infty,c}(\Omega)^n} \\ \grad|_{C_{\infty,c}(\Omega)} & 0 
 \end{pmatrix} \subseteqq \s A \subseteqq \begin{pmatrix}
    0 & \dive \\ {\grad} & 0 
 \end{pmatrix}.
\]
Now, $\s A$ is maximal monotone. Thus, $\s A$ is closed and we get that
\[
 \begin{pmatrix}
    0 & \interior\dive \\ \interior\grad & 0 
 \end{pmatrix} \subseteqq \s A \subseteqq \begin{pmatrix}
    0 & \dive \\ {\grad} & 0 
 \end{pmatrix}.
\]
The latter implies 
\[
 \iota_N^*\begin{pmatrix}
    0 & \interior\dive \\ \interior\grad & 0 
 \end{pmatrix}\iota_N \subseteqq \iota_N^* \s A \iota_N \subseteqq \iota_N^*\begin{pmatrix}
    0 & \dive \\ {\grad} & 0 
 \end{pmatrix}\iota_N.
\]
 From $\begin{pmatrix}
    0 & \interior\dive \\ \interior\grad & 0 
 \end{pmatrix}\subseteqq \s A$ it follows that \[N(\s A)^\bot \subseteqq N(\interior{\grad})^\bot \oplus N(\interior{\dive})^\bot.\] Thus, \[
     \iota_N^*\s A\iota_N \subseteqq \iota_N^*\begin{pmatrix}
    0 & \dive|_{N(\interior{\dive})^\bot} \\ {\grad}|_{N(\interior{\grad})^\bot} & 0 
 \end{pmatrix}\subseteqq \begin{pmatrix}
    0 & \dive|_{N(\interior{\dive})^\bot} \\ {\grad} & 0 
 \end{pmatrix}.\]
Now, let $(\phi,\psi)\in D(\iota_N\s A\iota_N^*)$. The latter inclusion shows that it suffices to estimate the norm of $\psi$ in the space $L_\nu^2(\R;H_1(\Omega))$. Moreover, we also read off that $\psi\bot N(\interior{\dive})$. Thus, $\psi\in R(\grad)$, which implies that $\psi$ takes almost everywhere values in the domain of $\curl$ and that $\curl \psi=0$. Recall 
\[
  a(\partial_0^{-1})\phi-\psi \in L_\nu^2(\R;D(\interior\dive)). 
\]
Usinge the smoothness assumptions on $a$, we compute
\begin{align*}
 \curl\left(a(\partial_0^{-1})\phi-\psi\right)& = \curl\left(a(\partial_0^{-1})\phi\right) \\
							     & = \curl\left(a(\partial_0^{-1})\right)\phi+a(\partial_0^{-1})\times \grad \phi \end{align*}
and
\[
  \dive (a(\partial_0^{-1})\phi)=\left(\dive a(\partial_0^{-1})\right)\phi+a(\partial_0^{-1})\grad \phi.
\]
Hence, with Theorem \ref{thm:gaffney}, we estimate pointwise almost everywhere
\begin{align*}
  &\abs{\psi}_{H_1}-\abs{a(\partial_0^{-1})\phi}_{H_1}\\&\leqq \abs{\psi-a(\partial_0^{-1})\phi}_{H_1}\\
 &\leqq c\left( \abs{\psi-a(\partial_0^{-1})\phi}_{L^2}+\abs{\dive\left(\psi-a(\partial_0^{-1})\phi\right)}_{L^2}+\abs{\curl\left(\psi-a(\partial_0^{-1})\phi\right)}_{L^2}\right)\\
 &\leqq c\left(\abs{\psi}_{L^2}+\abs{a(\partial_0^{-1})\phi}_{L^2}\right.\\ &\quad +\abs{\dive \psi}_{L^2}+\abs{\left(\dive a(\partial_0^{-1})\right)\phi}_{L^2}+\abs{a(\partial_0^{-1})\grad \phi}_{L^2}\\	&\quad \left. +\abs{\curl\left(a(\partial_0^{-1})\right)\phi}_{L^2}+\abs{a(\partial_0^{-1})\times \grad \phi}_{L^2}\right).  
\end{align*}
Thus, we get for some constant $c'>0$ that 
\begin{multline*}
  \abs{\psi}_{L^2_\nu(\R;H^1(\Omega))}\leqq c'\left(\abs{a(\partial_0^{-1})\phi}_{L_\nu^2(\R;H_1(\Omega))}+ \abs{\psi}_{L_\nu^2(\R;L^2(\Omega))}+\abs{a(\partial_0^{-1})\phi}_{L_\nu^2(\R;L^2(\Omega))}\right.\\  +\abs{\dive \psi}_{L_\nu^2(\R;L^2(\Omega))}+\abs{\left(\dive a(\partial_0^{-1})\right)\phi}_{L_\nu^2(\R;L^2(\Omega))}\\+\abs{a(\partial_0^{-1})\grad \phi}_{L_\nu^2(\R;L^2(\Omega))}+\abs{\curl\left(a(\partial_0^{-1})\right)\phi}_{L_\nu^2(\R;L^2(\Omega))}  \\ \left.	+\abs{a(\partial_0^{-1})\times \grad \phi}_{L_\nu^2(\R;L^2(\Omega))}\right). 
\end{multline*}
The smoothness assumptions on $a$ yield the assertion.
\end{proof}

Now, we are in the position to formulate the continuous dependence result. For simplicity, we only treat the case, where the operators in the material law do not depend on the spatial variables. The full homogenization problem will be discussed in future work. We adopt the strategy described in \cite[Section 1]{Waurick2012}. More specifically, we will treat the case of the particular class of operators being functions of $\partial_0^{-1}$ as discussed in Remark \ref{rem:mat_law_function}.

\begin{Sa}\label{thm:homog_impedance_higher}  Let $\nu>0, r>\frac{1}{2\nu}$. Assume that $\Omega\subseteqq \R^3$ is of class $C_5$ and such that $H_1(\Omega)\hookrightarrow\hookrightarrow L^2(\Omega)$. Let $\left(M_k\right)_k$ be a bounded sequence in $\s H^\infty(B(r,r))$ and denote $\s M_k\coloneqq M_k(\partial_0^{-1})$, $k\in\N$. If the conditions \eqref{eq:a_sksa} and \eqref{eq:AisMaxMon} hold and, in addition, there exists $c>0$ such that
\[
   \Re \langle z^{-1} M_k(z) u, u\rangle \geqq c\langle u,u\rangle \quad (k\in \N,\ z\in B(r,r),\ u\in L^2(\Omega)^4),
\]
then there is a subsequence $(n_k)_k$ of $(n)_n$ such that 
\[
   \partial_0 \s M_{n_k} + \s A \stackrel{G}{\to} \partial_0 \s M + \s A
\]
in $L_\nu^2(\R;L^2(\Omega)^4)$, where
\[
  \s M = \begin{pmatrix} \tau_{\textnormal{w}}\textnormal{-}\lim_{k\to\infty}\iota_N^* \s M_{n_k} \iota_N & 0 \\
                              0 & \partial_0^{-1}\left(\tau_{\textnormal{w}}\textnormal{-}\lim_{k\to\infty}\left(\partial_0\kappa_N^* \s M_{n_k} \kappa_N\right)^{-1}\right)^{-1} \\
              \end{pmatrix},
\]
with $\kappa_N\colon N(\s A)\to L_\nu^2(\R;L^2(\Omega)\oplus L^2(\Omega)^3)$ being the canonical injection. 
\end{Sa}
\begin{proof}
At first, we use Theorem \cite[Lemma 3.5]{PicPhy} to deduce that $\left(\partial_0\s M_k(\partial_0^{-1})\right)^{-1}$ is causal. Thus, from \cite[Lemma 3.8]{Trostorff2013}, we get that $\s M_k$ satisfies the positive definiteness condition imposed in Theorem \ref{thm:basic_conv_thm}, $k\in \N$.
 
 Let $f\in L_\nu^2(\R;L^2(\Omega)\oplus L^2(\Omega)^3)$ and consider the sequence $(u_n)_n$ in $L_\nu^2(\R;L^2(\Omega)\oplus L^2(\Omega)^3)$ satisfying
\[
   \left(\partial_0 \s M_n + \s A\right)u_n = f. 
\]
 The latter equation then reads as (note that being functions of $\s A$ the operators $\kappa_N, \iota_N$ commute with $\partial_0$):
\[
   \left(\partial_0 \begin{pmatrix} \iota_N^* \s M_n \iota_N & \iota_N^* \s M_n \kappa_N \\
                              \kappa_N^* \s M_n \iota_N & \kappa_N^* \s M_n \kappa_N \\
              \end{pmatrix}+
\begin{pmatrix} \iota_N^* \s A \iota_N & \iota_N^* \s A \kappa_N \\
                              \kappa_N^* \s A \iota_N & \kappa_N^* \s A \kappa_N \\
              \end{pmatrix}\right) \begin{pmatrix} \iota_N^* u_n\\ \kappa_N^* u_n \end{pmatrix} =\begin{pmatrix} \iota_N^* f\\ \kappa_N^* f \end{pmatrix}.
\]
Now, since $\partial_0^{-1}$ commutes with $\s A$, the $\s M_n$'s commute with $\kappa_N$ and $\iota_N$. Moreover, the skew-selfadjointness of $\s A$ implies that $\s A$ reduces $N(\s A)^\bot$. Thus, the latter system may be written as
\[
   \left(\partial_0 \begin{pmatrix} \iota_N^* \s M_n \iota_N & 0\\
                              0& \kappa_N^* \s M_n \kappa_N \\
              \end{pmatrix}+
\begin{pmatrix} \iota_N^* \s A \iota_N & 0\\
                              0& 0\\
              \end{pmatrix}\right) \begin{pmatrix} \iota_N^* u_n\\ \kappa_N^* u_n \end{pmatrix} =\begin{pmatrix} \iota_N^* f\\ \kappa_N^* f \end{pmatrix}.
\]
The latter gives the two (decoupled) equations:
\[
   \left(\partial_0\iota_N^* \s M_n \iota_N + \iota_N^* \s A \iota_N\right) \iota_N^* u_n = \iota_N^*f  
\]
and
\[
   \partial_0\kappa_N^* \s M_n \kappa_N \kappa_N^* u_n = \kappa_N^*f  
\] 
For the first equation, we use Theorem \ref{thm:basic_conv_thm}, the convergence of the equation in the stated manner follows from sequential compactness of bounded subsets of bounded linear operators in the weak operator topology.  
\end{proof}
\begin{rems}
 If, in the latter theorem, we restrict ourselves to the Hilbert space $N(\s A)^\bot$, i.e., using right-hand-sides, which are in $N(\s A)^\bot$, then the term involving $\kappa$ vanishes.  
\end{rems}

%
%
%

\subsection*{Applications to a singular perturbation problem}

To illustrate the applicability of Theorem \ref{thm:Solution_theory}, we give the following example of an elliptic/parabolic type equation, which is adopted from an example given in \cite{RainerPicard2013}. For this let $\Omega\subseteqq \R^n$ be open, bounded and connected and let $-\Delta$ be the Dirichlet Laplacian in $L^2(\Omega)$. Then $-\Delta$ is continuously invertible with compact resolvent. Let $\lambda \in (0,\lambda_1)$ for $\lambda_1$  being the smallest eigenvalue of $-\Delta$. Then, in particular, the operator $-\Delta-\lambda$ is maximal monotone. Now, let $\Omega_{\textnormal{p}},\Omega_{\textnormal{e}}\subseteqq \Omega$ be disjoint, measurable and such that $\Omega_{\textnormal{p}}\cup\Omega_{\textnormal{e}}=\Omega$. We let $\phi\colon \R\to\R$ be such that $\phi|_{(-\infty,0]}=0$, $\phi|_{(0,1)}=\id_{(0,1)}$ and $\phi|_{[1,\infty)}=1$. In $L_{\nu}^2(\R;L^2(\Omega))$, we consider for $\eps>0$ and given $f\in L_{\nu}^2(\R;L^2(\Omega))$ the problem of finding $u_\eps$ such that
\begin{equation}\label{eq:sing_pert}
   \left(\eps \partial_0 \phi(m_0)\1_{\Omega_\textnormal{p}}(m)+\1_{\Omega_\textnormal{e}}(m)(1-\phi(m_0))\tau_{-\eps} - \Delta\right) u_\eps =f,
\end{equation}
where $\tau_{-\eps}$ denotes the time-shift operator $\tau_{-\eps}g\coloneqq g(\cdot-\eps)$ for suitable $g$.
  At first, note that the latter problem is clearly well-posed. Indeed, this follows from 
\begin{multline*}
  \Re \langle \left(\eps \partial_0 \phi(m_0)\1_{\Omega_\textnormal{p}}(m)+\1_{\Omega_\textnormal{e}}(m)(1-\phi(m_0))\tau_{-\eps}+\lambda\right)u,\1_{\R_{\leqq a}}(m_0)u\rangle \\ \geqq \lambda' \langle u,\1_{\R_{\leqq a}}(m_0)u\rangle \quad (a\in \R)
\end{multline*}
for $u\in D(\partial_0)$, $\nu$ large enough and some $\lambda'\in (0,\lambda)$. On $\Omega_p$ the equation \eqref{eq:sing_pert} is of parabolic type and on $\Omega_e$ it is of elliptic type with an additional temporal variable. With 
\[
   \s M_\eps \coloneqq \eps \phi(m_0)\1_{\Omega_\textnormal{p}}(m)+\partial_0^{-1}\1_{\Omega_\textnormal{e}}(m)(1-\phi(m_0))\tau_{-\eps}+\partial_0^{-1}\lambda, 
\]
we get that $\s M_\eps \stackrel{\tau_{\textnormal{s}}}{\to} \s M_0=\partial_0^{-1}\1_{\Omega_\textnormal{e}}(m) (1-\phi(m_0))+\partial_0^{-1}\lambda$, where we denoted by $\tau_{\textnormal{s}}$ the strong operator topology. As strong convergence implies convergence in the weak operator topology, we infer with the help of Theorem \ref{thm:basic_conv_thm} that $(u_\eps)_{\eps>0}$ weakly converges as $\eps\to 0$ to the solution $u_0$ of the problem
\[
 \left(\1_{\Omega_\textnormal{e}}(m)(1-\phi(m_0)) - \Delta\right) u_0 =f,
\]
which itself is of pure elliptic type.

\section{Proof of Theorem \ref{thm:basic_conv_thm}}\label{sec:proof_of_main result}

For the proof we need several preparations.

\begin{Sa}[Theorem of Aubin-Lions, {\cite[p. 67, $2^\circ$]{Simon1987}}]\label{thm:Aubin-Lions} Let $H, K$ be Hilbert spaces, $I\subseteqq \R$ bounded, open interval. Assume that $K\hookrightarrow \hookrightarrow H$. Then 
\[
   H_1(I;H)\cap L^2(I;K)\hookrightarrow\hookrightarrow L^2(I;H).
\] 
\end{Sa}

\begin{Le}\label{le:weak_limits_and_causality} Let $H$ be a Hilbert space, $\nu>0$. Let $\left(\s M_n\right)_n$ be $\tau_{\textnormal{w}}$-convergent sequence in $L(L_\nu^2(\R;H))$ with limit $\s M$. If $\s M_n$ is causal for all $n\in\N$ then so is $\s M$. 
\end{Le}
\begin{proof}
 It suffices to observe that $\1_{\R_{<a}}(m_0)\in L(L_\nu^2(\R;H))$ for all $a\in \R$. Thus, the equation
\[
   \1_{\R_{<a}}(m_0)\s M_n  = \1_{\R_{<a}}(m_0)\s M_n \1_{\R_{<a}}(m_0)   
\]
carries over to the limit as $n\to\infty$.
\end{proof}

\begin{Sa}[weak-strong principle]\label{thm:wsp} Let $H,K$ be Hilbert spaces and such that $K\hookrightarrow\hookrightarrow H$. Let $\nu>0$ and $(v_n)_n$ be a weakly convergent sequence in $L^2_\nu(\R;K)\cap H_{\nu,1}(\R;H)$. Assume further that $\inf_{n\in\N} \inf \spt v_n>-\infty$. If $(\s M_n)_n$ is a $\tau_{\textnormal{w}}$-convergent sequence of causal operators in $L(L_\nu^2(\R;H))$, then
\[
   \textnormal{w-}\lim_{n\to\infty} \s M_n v_n = \left(\textnormal{$\tau_{\textnormal{w}}$-}\lim_{n\to\infty}\s M_n\right)\left(\textnormal{w-}\lim_{n\to\infty}v_n\right)\in L_\nu^2(\R;H).
\] 
\end{Sa}
\begin{proof} The uniform boundedness principle implies that both $(v_n)_n$ and $(\s M_n)_n$ are bounded sequences in $L^2_\nu(\R;K)\cap H_{\nu,1}(\R;H)$ and $L(L_\nu^2(\R;H))$, respectively. Thus, there exists a subsequence $(\s M_{n_k}v_{n_k})_k$ of $(\s M_{n}v_{n})_n$ which weakly converges to some $w\in L^2_\nu(\R;H)$. It suffices to identify $w$. For this, let $\phi \in C_{\infty,c}(\R;H)$ and define $a\coloneqq \sup \spt \phi$. Choose $\psi \in C_\infty(\R)$ such that $0\leqq \psi\leqq 1$,  $\psi = 1$ on $\R_{<a + 1}$ and $\psi=0$ on $\R_{>a+ 2}$. We denote $v\coloneqq \textnormal{w-}\lim_{n\to\infty}v_n$ and $\s M\coloneqq \textnormal{$\tau_{\textnormal{w}}$-}\lim_{n\to\infty}\s M_n$. Now, by Theorem \ref{thm:Aubin-Lions}, we deduce that $(\psi(m_0)v_n)_n$ converges to $\psi(m_0)v$ in $L_\nu^2(\R;H)$. For $n\in\N$, we compute
\begin{align*}
  \langle \s M_n v_n , \phi\rangle_{\nu,0} & = \langle \s M_n v_n, \1_{\R_{<a+1}}(m_0)\phi\rangle_{\nu,0} \\
                                   & = \langle \1_{\R_{<a+1}}(m_0) \s M_n v_n, \phi\rangle_{\nu,0} \\
                                   & = \langle \1_{\R_{<a+1}}(m_0) \s M_n \1_{\R_{<a+1}}(m_0) v_n, \phi\rangle_{\nu,0} \\
                                   & = \langle \1_{\R_{<a+1}}(m_0) \s M_n \1_{\R_{<a+1}}(m_0)\psi(m_0) v_n, \phi\rangle_{\nu,0} \\
                                   & = \langle \1_{\R_{<a+1}}(m_0) \s M_n \psi(m_0) v_n, \phi\rangle_{\nu,0} \\
                                   & = \langle \s M_n \psi(m_0) v_n,\1_{\R_{<a+1}}(m_0)  \phi\rangle_{\nu,0} \\
                                   & \to \langle \s M \psi(m_0) v,\1_{\R_{<a+1}}(m_0)  \phi\rangle_{\nu,0} \\
                                   & = \langle \s M v, \phi\rangle_{\nu,0},
\end{align*}where we have used that the $\s M_n$'s and $\s M$ are causal, see also Lemma \ref{le:weak_limits_and_causality}. Hence, 
\[
   \langle w,\phi\rangle_{\nu,0} = \langle \s M v,\phi\rangle_{\nu,0}
\]
for all $\phi\in C_{\infty,c}(\R;H)$. Thus, $w=\s M v$.
\end{proof}

\begin{rem}
  The support condition for the $v_n$'s is needed to make Theorem \ref{thm:Aubin-Lions} applicable.
\end{rem}

\begin{Le}\label{le:weak_limits_and_commutator} Let $H$ be a Hilbert space, ${\s D}$ densely defined, closed, linear operator in $H$ with $0\in\rho({\s D})$. Let $(M_n)_n$ be a sequence in $L(H)$ converging in the weak operator topology to some $M$ and having bounded commutators with $\s D$. Then
\[
   [M_n,{\s D}]\to \overline{M\s D-\s DM} \quad(n\to\infty)
\]
 in the weak operator topology. In particular, $M\s D-\s DM$ extends to a bounded linear operator, $M$ has a bounded commutator with $\s D$ and
\[
   \s D M_n u \rightharpoonup \s D M u \quad(n\to\infty,\ u\in D(\s D)).
\]
\end{Le}
\begin{proof}
 For $x,y\in H$, $n\in\N$ we compute
\begin{align*}
  \langle [M_n,{\s D}]{\s D}^{-1}x,({\s D}^{-1})^*y \rangle &=\langle (M_n{\s D}-{\s D}M_n){\s D}^{-1}x,({\s D}^{-1})^*y\rangle \\
                                             &= \langle {\s D}^{-1}(M_n{\s D}-{\s D}M_n){\s D}^{-1}x,y\rangle \\
                                             &= \langle ({\s D}^{-1}M_n-M_n{\s D}^{-1})x,y\rangle \\
                                             &= \langle {\s D}^{-1}M_nx,y\rangle -\langle M_n{\s D}^{-1}x,y\rangle\\
                                             &\to \langle {\s D}^{-1}Mx,y\rangle -\langle M{\s D}^{-1}x,y\rangle\\
                                             &= \langle [M,{\s D}]{\s D}^{-1}x,({\s D}^{-1})^*y \rangle.  
\end{align*}
By the boundedness of $([M_n,{\s D}])_n$ and the density of both $D({\s D})$ and $D({\s D}^*)$, we get the first convergence result. In order to see the last convergence result, we compute for $n\in\N$ and $u\in \s D$:
\[
   \s D M_n u=[\s D,M_n]u+M_n\s D u \rightharpoonup [\s D,M]u+M\s D u =\s D M u \qedhere 
\]
\end{proof}

\begin{Le}\label{le:commutator_is_causal} Let $H$ be Hilbert space, $\nu>0$. Let $M\in L(L_\nu^2(\R;H))$ be causal and such that $M$ has a bounded commutator with $\partial_0$. Then $[M,\partial_0]$ is causal. 
\end{Le}
\begin{proof}
  Let $a\in \R$ and $\phi\in C_{\infty,c}(\R;H)$ be such that $\spt \phi\geqq a$. Thus, $\spt M\phi\geqq a$ and $\spt \phi'\geqq a$. Now, since $M\phi\in D(\partial_0)$ since $M[D(\partial_0)]\subseteqq D(\partial_0)$ we further get $\spt \partial_0 M\phi\geqq a$. Hence, we arrive at $\spt [M,\partial_0]\phi= \spt \left(M\partial_0-\partial_0M\right)\phi\geqq a$. The continuity of $[M,\partial_0]$ together with Remark \ref{rem:causality} imply the assertion.
\end{proof}

\begin{Fo}\label{cor:wsp_appl} Let $K,H$ be Hilbert spaces, $\nu>0$. Let $\left(\s M_n\right)_n$ be a bounded sequence of causal mappings in $L(L_\nu^2(\R;H))$ converging in the weak operator topology to some $\s M$ and having bounded commutators with $\partial_0$. Assume that $K\hookrightarrow\hookrightarrow H$. Let $(u_n)_n$ be a weakly convergent sequence in $H_{\nu,-1}(\R;K)\cap L_\nu^2(\R;H)$ with limit $u$ and such that $\inf_{n\in\N}\inf u_n>-\infty$.
Then $\s M_n u_n \rightharpoonup \s M u$ in $L_\nu^2(\R;H)$ as $n\to\infty$. 
\end{Fo}
\begin{proof}
 At first note that $(\partial_0^{-1}u_n)_n$ is weakly convergent in $L^2_{\nu}(\R;K)\cap H_{\nu,1}(\R;H)$. Moreover, $[\s M_n,\partial_0]$ is causal by Lemma \ref{le:commutator_is_causal} for all $n\in\N$. Furthermore, $[\s M_n,\partial_0]\stackrel{\tau_{\text{w}}}{\to} [\s M,\partial_0]$, by Lemma \ref{le:weak_limits_and_commutator}. Thus, for $n\in\N$ we deduce with the help of Theorem \ref{thm:wsp} that
\begin{multline*}
  \s M_n u_n = \s M_n \partial_0 \partial_0^{-1}u_n = [\s M_n,\partial_0]\partial_0^{-1}u_n+\partial_0 \s M_n \partial_0^{-1}u_n\\
 \rightharpoonup [\s M,\partial_0]\partial_0^{-1}u + \partial_0 \s M\partial_0^{-1} u =\s M u\in L^2_{\nu}(\R;H),
\end{multline*}
where we have used that $\partial_0^{-1}$ is weakly continuous and causal.
\end{proof}

\begin{Le}\label{le:convergence_of_unbounded_A} Let $H$ be a Hilbert space, $\nu>0$. Let $\s A$ be a densely defined, closed, linear operator in $L_\nu^2(\R;H)$ with $0\in \rho(\s A)$. Assume that $\partial_0^{-1}\s A^{-1}=\s A^{-1}\partial_0^{-1}$.  Let $(u_n)_n$ be a bounded sequence in $H_{-1,1}(\partial_0,\s A)\cap L_\nu^2(\R;H)$, which weakly converges in $L_\nu^2(\R;H)$ to $u\in L_\nu^2(\R;H)$. Then $u\in H_{-1,1}(\partial_0,\s A)$ and
\[
   \s A u_n \rightharpoonup \s A u \in H_{-1,0}(\partial_0,\s A).
\]
\end{Le}
\begin{proof}
 Let $(u_{n_k})_k$ be a weakly convergent subsequence of $(u_n)_n$ in $H_{-1,1}(\partial_0,\s A)$. Denote its limit by $w$. Note that $\partial_0^{-1}u_n \rightharpoonup \partial_0^{-1}u \in H_{\nu,1}(\R;H)\hookrightarrow L_\nu^2(\R;H)$, by unitarity of $\partial_0^{-1}$. Moreover, by Remark \ref{rem:SoLa}, $\partial_0^{-1}\colon H_{-1,1}(\partial_0,\s A)\to H_{0,1}(\partial_0,\s A)$ is unitary. Hence, we get that $\partial_0^{-1}u_{n_k}\rightharpoonup \partial_0^{-1}w\in H_{0,1}(\partial_0,\s A)\hookrightarrow L_\nu^2(\R;H)$. Hence, $\partial_0^{-1}w=\partial_0^{-1}u$. Thus, $u=w$ and $(u_n)_n$ weakly converges in $H_{-1,1}(\partial,\s A)$. Now, by Remark \ref{rem:SoLa}, the operator $\s A\colon H_{-1,1}(\partial_0,\s A)\to H_{-1,0}(\partial_0,\s A)$ is continuous. Thus, we deduce the asserted convergence. 
\end{proof}

\begin{proof}[Proof of Theorem \ref{thm:basic_conv_thm}] The well-posedness of the limiting equation, i.e., continuous invertibility and causality of $(\partial_0\s M+\s A)$ in $L_\nu^2(\R;H)$ follows from Lemma \ref{le:weak_limits_and_commutator} together with Theorem \ref{thm:Solution_theory}.

Now, we prove the version, which is asserted in Remark \ref{rem:stronger_conv}. Let $(f_n)_n$ in $L_\nu^2(\R;H)$ be a weakly convergent sequence with $\inf_n\inf \spt f_n>-\infty$; we denote its limit by $f$. For $n\in\N$ we define
\[
   u_n\coloneqq \left(\partial_0 \s M_n +\s A\right)^{-1}f_n.
\]
By causality (see Theorem \ref{thm:Solution_theory} and Remark \ref{rem:causality}), we get that \[\inf_{n\in\N}\inf \spt u_n\geqq \inf_n\inf\spt f_n>-\infty.\] Moreover, $(u_n)_n$ is bounded in $L_\nu^2(\R;H)\cap H_{-1,1}(\partial_0,\s A+1)$ by Corollary \ref{cor:cont_estim} and the uniform boundedness principle applied to $(f_n)_n$. Now, let $(u_{n_k})_k$ be a $L_\nu^2(\R;H)$-weakly convergent subsequence of $(u_n)_n$. We denote the respective limit by $u$. Now, for $n\in\N$ we have
\begin{equation}\label{eq:nk}
   \partial_0 \s M_{n_k} u_{n_k} + \s A u_{n_k} = f_{n_k}
\end{equation}
in $H_{-1,0}(\partial_0,\s A+1)$. Now, by Corollary \ref{cor:wsp_appl} for the first term and Lemma \ref{le:convergence_of_unbounded_A} for the second term on the left side of equation \eqref{eq:nk}, we may let $k\to\infty$ in \eqref{eq:nk}. We arrive at
\[
   \partial_0 \s M u + \s A u = f
\]
 in $H_{-1,0}(\R;H)$. Moreover, by construction, $u\in L_\nu^2(\R;H)$ and $(\partial_0\s M+\s A)u=f\in H_{0,0}(\partial_0,\s A+1)$. Thus, $u\in D(\partial_0\s M+\s A)$, by Remark \ref{rem:domainoftime-space}. Now, since $\left(\partial_0\s M + \s A\right)$ is continuously invertible in $L_\nu^2(\R;H)$ the sequence $(u_n)_n$ itself weakly converges. 

In order to see that Theorem \ref{thm:basic_conv_thm} holds, apply the previous part to constant sequences $(f_n)_n=(f)_n$ for some $f\in C_{\infty,c}(\R;H)$. It remains to observe that $C_{\infty,c}(\R;H)$ is dense in $L_\nu^2(\R;H)$ and that $\left(\left(\partial_0 \s M_n +\s A\right)^{-1}\right)_n$ is bounded. 
\end{proof}

\section*{Acknowledgements}

The author is indebted to Rainer Picard and Sascha Trostorff for fruitful discussions particularly concerning the example of the wave equation with impedance type boundary conditions and for useful comments to streamline the present exposition.


\end{document}